\documentclass[12pt,letterpaper,reqno]{amsart}
\usepackage{graphicx}
\usepackage{url}
\usepackage{ bbold }
\usepackage{hyperref}
\RequirePackage[margin=1in]{geometry}
\RequirePackage[english]{babel}
\RequirePackage[utf8]{inputenc}
\RequirePackage{mathtools}
\RequirePackage{amsmath}
\RequirePackage{amssymb}
\RequirePackage{amsfonts}
\RequirePackage{amsthm}
\RequirePackage{mathrsfs}
\RequirePackage{cite}
\RequirePackage{url}
\RequirePackage{hyperref}
\RequirePackage[labelfont=bf,labelsep=period,justification=raggedright]{caption}
\RequirePackage{subcaption}
\RequirePackage{indentfirst}
\RequirePackage[colorinlistoftodos]{todonotes}
\RequirePackage{tikz}
\RequirePackage{tkz-fct}
\RequirePackage{tikz-3dplot}
\RequirePackage[all]{xy}
\RequirePackage{float}
\RequirePackage{placeins}
\RequirePackage{mdframed}
\RequirePackage{multicol}
\RequirePackage{pifont}
\RequirePackage{lipsum}
\RequirePackage{algorithm}
\RequirePackage{algpseudocode}
\usetikzlibrary{calc}
\allowdisplaybreaks

\renewcommand{\lg}[2]{\left(\frac{#1}{#2}\right)}

\theoremstyle{plain}
\newtheorem{thm}{Theorem}
\newtheorem*{thm*}{Theorem}
\newtheorem{lemma}[thm]{Lemma}
\newtheorem*{lemma*}{Lemma}

\newtheorem*{cor*}{Corollary}
\newtheorem{conj}[thm]{Conjecture}
\newtheorem*{conj*}{Conjecture}
\newtheorem{prop}[thm]{Proposition}
\newtheorem*{prop*}{Proposition}

\newtheorem*{qn*}{Question}

\newtheorem*{claim*}{Claim}
\newtheorem{prob}{Problem}
\newtheorem*{prob*}{Problem}

\theoremstyle{definition}
\newtheorem{defn}[thm]{Definition}
\newtheorem*{defn*}{Definition}

\newtheorem*{sol*}{Solution}

\newtheorem*{axm*}{Axiom}

\theoremstyle{remark}
\newtheorem{rem}[thm]{Remark}
\newtheorem*{rem*}{Remark}

\newenvironment{fthm*}
    {\begin{mdframed}[linewidth=2.0pt]\begin{thm*}}
    {\end{thm*}\end{mdframed}}

\newenvironment{fdefn*}
    {\begin{mdframed}[linewidth=2.0pt]\begin{defn*}}
    {\end{defn*}\end{mdframed}}

\newenvironment{fprob*}
    {\begin{mdframed}[linewidth=2.0pt]\begin{prob*}}
    {\end{prob*}\end{mdframed}}

\renewcommand{\pagebreak}{\clearpage}

\numberwithin{equation}{section}
\numberwithin{thm}{section}

\title{Murmurations in Short Intervals for Elliptic Curves Ordered by Height}

\author[Aarya J. Kumar]{Aarya J. Kumar}
\address[Aarya J. Kumar]{Princeton University}
\email{ajfromny2004@gmail.com}

\begin{document}

\begin{abstract}
    Murmurations are the phenomenon, first observed by He, Lee, Oliver, and Pozdnyakov (2025), that the average value of $\epsilon(E) a_p(E)$, the $p$-th $L$-function coefficient weighted by the root number, over elliptic curves of conductor in some range $N(E) \in [N_1, N_2]$ seem to approximate a continuous, oscillating function. While murmurations have been proven for other families of $L$-functions (notably by Zubrilina (2025) for weight $k$ modular forms), the only unconditional murmuration formula for elliptic curves is a partial result by Sawin and Sutherland (2025) ordering by naive height instead of conductor. In this setting, murmurations appear by first averaging $\epsilon(E) a_p(E)$ for a single curve over $p \in (C_1, C_2)$ and then average over elliptic curves of bounded height. However, Sawin and Sutherland's results average $\epsilon(E) a_n(E)$ over all of $\mathbb{N}$ instead of a small interval $(C_1, C_2).$ In this work, we prove how small the interval $(C_1, C_2)$ can be taken while preserving the conclusion of Sawin and Sutherland's first murmuration formula (Proposition 7).
\end{abstract}

\maketitle

\section{Introduction}

Murmurations are a statistical phenomenon, first observed by He, Lee, Oliver, and Pozdnyakov for elliptic curves in \cite{HLOP25}. Specifically, they concern oscillations in the averages of $L$-function coefficients $a_p$ (for fixed prime $p$) over curves with similar conductor. These oscillations, which are dependent on the rank of the averaged curves, appear to approximate a continuous function. Recent work has demonstrated that the phenomenon holds more generally for $L$-functions of families of modular forms. 

Results have been proven regarding murmurations of other families of $L$-functions, in particular Zubrilina's murmuration density formula \cite{Zub25} for modular forms of even weight $k$. However, no explicit proof of a murmuration density formula has been proved for $L$-functions of elliptic curves. The work of Sawin and Sutherland, specifically \cite[Theorem 2]{SS25}, marks the first unconditional result in the direction of such a murmuration density formula. Based on this theorem, Sawin and Sutherland form a hypothesis \cite[Conjecture 1]{SS25} for what the murmuration density formula should be, and they prove that their given formula is indeed a continuous function. We remark that in order to enable the computation of murmurations of elliptic curves, the formulation of the problem in \cite{SS25} deviates from the original phenomenon in several ways, one being the choice to order elliptic curves by naive height rather than conductor. This change is significant, as averages of $a_p(E)$ over curves of bounded height no longer form visible murmuration patterns, and it becomes necessary instead to first average $a_p(E)$ for fixed curve $E$ over a range $p \in (C_1, C_2)$ and then average over elliptic curves $E.$

The goal of this paper is to prove Theorem \ref{Main Result}, which is a refinement of prior work in \cite{SS25} in the direction of \cite[Conjecture 1]{SS25}. In this introduction, we begin by recalling prior work concerning the phenomena of murmurations, primarily following \cite{SS25} since we aim to orient the discussion towards results in this paper that we will adapt. Finally, after discussing in detail the main results of \cite{SS25}, we state our main result, Theorem \ref{Main Result}, and provide an overview of the proof technique along with directions for further work.

Finally, a note on notation. Throughout the discussion, $E$ will refer to an elliptic curve over $\mathbb{Q}$ and $E_{a, b}$ will refer specifically to the short Weierstrass equation $y^2 = x^3 + ax + b.$ $N(E)$ will be used to represent the conductor of $E,$ and $H(E)$ will represent the naive height. If $S$ is a given family of elliptic curves and $f$ is a function or property of elliptic curves, we notate an average over $S$ as
$$\mathbb{E}_{E \in S} [f(E)] := \frac{1}{|S|} \sum_{E \in S} f(E).$$

\subsection{Murmurations of Elliptic Curves}

The discovery of murmurations of elliptic curves has its origins in attempts such as \cite{HLO20} to use machine-learning to predict the algebraic properties of an elliptic curve $E$, such as rank, from a finite number of its $L$-function coefficients $a_n(E)$. This algorithm was surprisingly efficient, prompting exploration in \cite{HLOP25} regarding whether the behavior of $L$-function coefficients is influenced by the ranks of elliptic curves. 

\begin{figure}[h]
    \centering
    \includegraphics[width=0.5\linewidth]{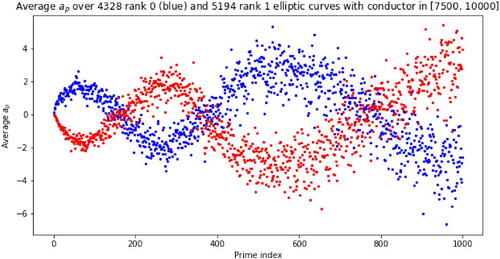}
    \caption{The blue plot records the average of $a_p$ over $4328$ random elliptic curves $E$ with $\text{rank}(E) = 0$ and $N(E) \in [7500, 10000],$ with $p \in [1, 1000].$ The red plot is the average of $a_p$ over $5194$ random curves $E$ with $\text{rank}(E) = 1$ and $N(E) \in [7500, 10000]$. Figure taken from \cite[Figure 6]{HLOP25}.}
    \label{fig:original murmuration example}
\end{figure}

Since, given a fixed $x$, the equation $y^2 = x^3 + ax + b$ over $\mathbb{F}_p$, has probability $\frac{1}{2}$ of yielding a solution for $y,$ the expected number of solutions for the curve $E_{a, b}$ over $\mathbb{F}_p$ is $p + 1$, including a point at infinity. One thus expects that the coefficients $a_n(E)$ should have an average value of $0$. However, \cite{HLOP25} demonstrates through empirical evidence that given a fixed rank $r \in \{0, 1\}$ and a range $[N_1, N_2]$ of natural numbers, the average
$$\mathbb{E}_{\{E: \text{rank}(E) = r, N(E) \in [N_1, N_2]\}}[a_p]$$
is not zero but instead exhibits oscillatory behavior dependent on $p.$ Moreover, it appears to be true that
$$\mathbb{E}_{\{E: \text{rank}(E) = 0, N(E) \in [N_1, N_2]\}}[a_p] \approx - \mathbb{E}_{\{E: \text{rank}(E) = 1, N(E) \in [N_1, N_2]\}}[a_p].$$
For an example of these phenomena, see Figure \ref{fig:original murmuration example}.

Further observations on murmurations were subsequently made by Sutherland in a letter \cite{Sut22} to Rubenstein and Sarnak. These include a so-called ``scale-invariance" property wherein the crossing points of the oscillations seems to be proportional to the interval $[N_1, N_2].$ Additionally, the ``sign" of these murmurations ought to depend not on the rank $r$ of the averaged set of elliptic curves but only on the parity $(-1)^{r},$ which, according to the Birch and Swinnerton-Dyer conjecture \cite{BSD65}, is equal to the root number $\epsilon(E)$ appearing in the functional equation
$$(2 \pi)^{-s} \Gamma(s) L(E, s) = \epsilon(E) N(E)^{1 - s} (2\pi)^{s - 2} \Gamma(2 - s)L(E, 2 - s).$$
The oscillations are not quite the same because curves of larger rank are expected to have larger values of $a_p(E).$ Nevertheless, these curves of high rank are sparse, expected to exist with zero density. Accordingly, rather than restricting to a single rank, Sutherland focuses on averages of the form $\mathbb{E}_{E: \epsilon(e) = r, N(E) \in [C_1, C_2]}[a_p]$ with $r \in \{\pm 1\}.$ In \cite{SS25}, Sawin and Sutherland consider the average
$$\mathbb{E}_{E: N(E) \in [C_1, C_2]}[\epsilon(E) a_p],$$
which weights elliptic curves by their root number. According to empirical data, these averages have the advantage of being smoother and converging more rapidly.

Perhaps most importantly, Sutherland observes that murmurations occur not only for elliptic curves but more generally for any family of primitive, normalized arithmetic $L$-functions, including the $L$-functions of modular forms. Recall that the $L$-function of a modular form $f(z) = \sum_{n = 0}^\infty a_n q^n$ (where $q = e^{2\pi i}$) is defined to be the Dirichlet series
$$L(f, s) := \sum_{n = 1}^\infty \frac{a_n}{n^{-s}}.$$
Notably, the case of $L$-functions $L = L(E, s)$ falls under the case of $L$-functions of modular forms by the modularity theorem, first proved by Wiles and Taylor in \cite{Wil95} and \cite{WT95} in the special case of semistable elliptic curves and then more generally by Breuil, Conrad, Diamond, and Taylor in \cite{BCDT01}.

Some work has been done to prove murmuration density formulae for various families of modular forms and automorphic forms. These cases turn out to be more manageable because trace formulae of automorphic forms are compatible with root numbers $\epsilon,$ so averages over $\epsilon^{-1} a_p$ turn out to be convenient. Numerous results have been established for families of modular forms of fixed weight and varying level, specifically Zubrilina's murmuration formula \cite{Zub25} for $L$-functions of modular forms of weight $k.$ Of these results, the case of $k = 2$ comes closest to providing a result about elliptic curves, since the $L$-functions of rational elliptic curves are a very strict subset of this family of modular $L$-functions by the modularity theorem. Cowan has proved a conjectural murmuration formula for elliptic curves in \cite{Cow25}, based on random matrix theory and contingent on the Ratios Conjectures \cite{DHP15}, the generalized Riemann hypothesis for various families of $L$-functions, and several other conjectures.

The work of Sawin and Sutherland \cite{SS25} marks the first attempt to establish an unconditional murmuration formula for elliptic curves. Since it is difficult to establish any statistical results for elliptic curves ordered by conductor, Sawin and Sutherland work with elliptic curves ordered by naive height
$$H(E_{A, B}) := \max(4|A|^3, 27|B|^2).$$
Previously, Sutherland \cite{Sut22} noticed that the existence of murmurations is specific to ordering by conductors and that oscillations are not visible in alternative orderings. However, \cite{SS25} establishes through empirical evidence that this issue can be amended through ``local averaging." That is, instead of averaging the coefficient $a_n(E)$ over a family of elliptic curves, we average the sum $\sum_{p \in (C_1 N(E), C_2 N(E)} [\epsilon(E) a_p(E)].$ As discussed in \cite{SS25}, it has been hypothesized by Jonathan Bober that the relationship between $a_p(E)$ and $\epsilon(E)$ seems to depend on the ratio $\frac{p}{N},$ further justifying the use of a sum of coefficients over such a range. The definition of murmurations that we will work from, following \cite[Equation 1]{SS25}, is as follows.

\begin{defn} \label{Def Murmurations}
    Murmurations for elliptic curves ordered by naive height are given by
    $$\mathbb{E}_{\{E : H(E) \leq X \}} \left[ \frac{\log\left( N(E) \frac{C_1 + C_2}{2}  \right)}{N(E)}  \sum_{p \in (C_1 N(E), C_2N(E) )} \epsilon(E) a_p(E) \right].$$
\end{defn}
The factor of $\frac{\log\left( N(E) \frac{C_1 + C_2}{2}  \right)}{N(E)}$ included in the above average is the reciprocal for the expected number of primes in the interval $(C_1 N(E), C_2 N(E))$ according to the prime number theorem. Sawin and Sutherland's formula for murmurations of elliptic curves is the following.
\begin{conj}[{\cite[Conjecture 1]{SS25}}]\label{Murmuration Formula Conjecture}
    For real numbers $0 < C_1 < C_2$ we have
    $$\lim_{X \rightarrow \infty} \mathbb{E}_{\{E : H(E) \leq X \}} \left[ \frac{\log\left( N(E) \frac{C_1 + C_2}{2}  \right)}{N(E)}  \sum_{p \in (C_1 N(E), C_2N(E) )} \epsilon(E) a_p(E) \right]$$
    $$= \int_{C_1}^{C_2} 2\pi \sqrt{u} \sum_{\substack{q \in \mathbb{Z}\\ \text{$q$ square-free}}} \sum_{m \in \mathbb{N}} \frac{\mu(\gcd(m, q))}{q m \phi \left( \frac{q}{\gcd(m, q)}  \right)} J_1\left( 4\pi \frac{\sqrt{u} m}{q}  \right) \prod_{p \mid q} \hat{\ell}_{p, 2v_p(m)} \prod_{p \mid m, p \nmid q} \ell_{p, 2v_p(m)} du.$$
\end{conj}
Sawin and Sutherland note that the right-hand side of the above equation is indeed a continuous function in the real variables $C_1$ and $C_2.$ In the above, $J_1$ is a Bessel function of the first kind, and $\ell_{p, v}$ and $\hat{\ell}_{p, v}$ are local constant factors dependent on primes $p$ and valuations $v.$ The following definition for $\ell_{p, v}$ is taken from \cite[Definition 5]{SS25}. The similar but more complicated version $\hat{\ell}_{p, v}$ is irrelevant for our purposes but can be found in \cite[Definition 6]{SS25}.

\begin{defn}[{\cite[Definition 5]{SS25}}]\label{local factor}
    The local factor $\ell_{p, v}$ for prime $p$ and valuation $v \in \mathbb{N} \cup\{0\}$ is defined as the integral
    $$\ell_{p, v} := \frac{1}{1 - p^{-10}} \int_{(A, B) \in \mathbb{Z}_p^2 \setminus(p^4 \mathbb{Z}_p \times p^6 \mathbb{Z}_p)} a_{p^v} (E_{A, B}),$$
    where $\mathbb{Z}_p$ denotes the $p$- adic integers and the integral is taken against the uniform measure on $\mathbb{Z}_p^2$ with total mass $1.$
\end{defn}
Intuitively, this quantity can be understood as a weighted average of all the possible $p^v$-th coefficients $a_{p^v}$. We will see later through Remark \ref{p adic coefficients} and the proof of Lemma \ref{Average of coefficients} that this integral really can be written as a finite average. Most of the evidence for \cite[Conjecture 1]{SS25} comes from the following theorem.

\begin{thm}[{\cite[Theorem 2.2]{SS25}}]\label{Murmurations Theorem 2}
    Let $W$ be a smooth, compactly supported function on $(0, \infty).$ Then the limit
    $$\lim_{P \rightarrow \infty} \lim_{X \rightarrow \infty} \mathbb{E}_{\{ E : H(E) \leq X\}} \left[ \frac{\prod_{p \leq P} \left( 1 - \frac{1}{p} \right)^{-1}}{N(E)} \sum_{\substack{n \in \mathbb{N} \\ \text{$p \nmid n$ for $p \leq P$}}} W \left( \frac{n}{N(E)}\right) \epsilon(E) a_n(E)  \right]$$
    exists and is equal to
    $$\int_0^\infty W(u) 2\pi \sqrt{u} \sum_{\substack{q \in \mathbb{N} \\ \text{$q$ square-free}}} \sum_{m \in \mathbb{N}} \frac{\mu(\gcd(m, q))}{q m \phi\left( \frac{q}{\gcd(m, q)}  \right)} J_1\left( 4\pi \frac{\sqrt{u} m}{q} \right) \prod_{p \mid q} \hat{\ell}_{p, 2v_p(m)} \prod_{p \mid m, p \nmid q} \ell_{p, 2v_p(m)} du.$$
\end{thm}

The result \cite[Theorem 2.2]{SS25} differs from \cite[Conjecture 1]{SS25} in several ways. The first is that rather than including a sum of $a_p(E)$ for all primes, \cite[Theorem 2.2]{SS25} sieves out finitely many primes at a time and then takes the limit as the number of primes sieved approaches infinity. While this sieve process is heuristically equivalent to summing over primes $a_p$, it is unlikely that an equivalence can be proved without extreme difficulty. 

Secondly, instead of summing over coefficients in a small interval $(C_1 N(E), C_2 N(E)),$ \cite[Theorem 2.2]{SS25} sums coefficients over all of $\mathbb{N}$ against a weight function of compact support. Here, one should imagine taking the weight function to be the indicator function of an interval, or more accurately, a smooth approximation of such an indicator function. However, \cite[Theorem 2.2]{SS25} does not adequately account for dependence on the interval $(C_1 N(E), C_2 N(E)).$ The problem of adjusting the results of \cite{SS25} to account for this issue is a more achievable problem and constitutes the focus of this paper.

\subsection{Main Results}

The remainder of this paper will be focused on adapting the following result of \cite{SS25}, which provides an analog for \cite[Theorem 2.2]{SS25} that sums over all coefficients $a_n$, not just coefficients of prime index.

\begin{prop}[{\cite[Proposition 7]{SS25}}] \label{Proposition 7}
    For $W$ a smooth, compactly supported function on $(0, \infty),$ the limit
    $$\lim_{X \rightarrow \infty} \mathbb{E}_{E: H(E) \leq X} \left[ \frac{1}{N(E)} \sum_{n = 1}^\infty W\left(   \frac{n}{N(E)} \right) \epsilon(E) a_n(E)\right]$$
    exists and evaluates to 
    $$\int_0^\infty W(u) 2\pi \sqrt{u} \sum_{m = 1}^\infty \frac{1}{m} J_1\left(4\pi \sqrt{u}m\right) \prod_{p \mid m} \ell_{p, 2v_p(m)} du.$$
\end{prop}
The problem of replacing the infinite sum over $n$ in \cite[Proposition 7]{SS25} with a sum over an interval of finite length is equivalent to the problem of replacing the test function $W$ with the test function $W(L(u - x)),$ which contracts the support of $W$ by a factor of $L$ and then shifts the function to the right by $Lx.$ Our primary challenge is to determine how large $L$ can be chosen relative to the bound $X$ on naive height. Since murmurations only occur over intervals of $a_n$ and not single values of $a_n,$ it is reasonable to ask at what point the relationship in \cite[Proposition 7]{SS25} breaks down. Our main result is Theorem \ref{Main Result}, which modifies \cite[Proposition 7]{SS25} in such a way.

\begin{thm} \label{Main Result}
    Let $W$ be a smooth, compactly supported function on $(0, \infty).$ Let $L = L(X)$ be such that $L \ll X^{\frac{3}{10} - \epsilon}$ for some $\epsilon > 0,$ and let $x$ be a fixed constant. Then the limit
    $$\lim_{X \rightarrow \infty} \mathbb{E}_{E: H(E) \leq X} \left[ \frac{1}{N(E)} \sum_{n = 1}^\infty W\left( L\left(  \frac{n}{N(E)} - x\right) \right) \epsilon(E) a_n(E)\right]$$
    exists and evaluates to
    \begin{equation}\label{RHS-of-main-theorem}
        \int_0^\infty W(L(u - x)) 2\pi \sqrt{u} \sum_{m = 1}^\infty \frac{1}{m} J_1\left(4\pi \sqrt{u}m\right) \prod_{p \mid m} \ell_{p, 2v_p(m)} du.
    \end{equation}
\end{thm}

The proof technique for Theorem \ref{Main Result} is similar to that of \cite[Proposition 7]{SS25}, which for the most part follows from an application of Voronoi summation. Most complications come from a single step in this proof, where one computes the average (for fixed $n \in \mathbb{N}$)
$$\lim_{X \rightarrow \infty} \mathbb{E}_{\{E: H(E) \leq X\}}[a_n(E)] = \prod_{p \mid n} \ell_{p, v_p(n)}.$$
In determining how large we can permit $L$ to be for the conclusion of Theorem \ref{Main Result} to hold, it is necessary to understand how large the error terms of this equality are. This question is answered by the following lemma.

\begin{lemma}\label{Average of coefficients}
    For $n \in  \mathbb{N},$
    \begin{align*}
        \mathbb{E}_{\{E: H(E) \leq X\}}[a_n(E)] &= \prod_{p \mid n} \ell_{p, v_p(n)} + O \left( X^{-\frac{2}{3}} n^{\frac{1}{2}} \right)\\
    &+ O_\epsilon\left(X^{-\frac{1}{3}} q^{\frac{1}{2} + \epsilon}\prod_{p : p^2 \mid n} p^{\frac{v_p(n)}{2} + \epsilon} \right) + O_\epsilon \left( X^{-\frac{3}{4}} q^{\frac{3}{2} + \epsilon}\prod_{p : p^2 \mid n} p^{\frac{v_p(n)}{2} + \epsilon} \right).
    \end{align*} 
\end{lemma}

Here, $q = q(n)$ is a constant dependent only on $n,$ defined in Definition \ref{def q}. It grows similarly to the square-free part of $n.$ If $L$ is large enough compared to $X,$ the contribution of $L$ in \ref{RHS-of-main-theorem} can end up being large enough that it overwhelms at least one of the two error terms in Lemma \ref{Average of coefficients}. In proving Lemma \ref{Average of coefficients}, we rely on standard analytic techniques such as the ``Polya-Vinogradov" method in Lemma \ref{Polya Vinogradov} to rewrite $\sum_{E: H(E) \leq X} a_n(E)$ as a sum of exponential terms, which combines error terms efficiently into geometric series.

Of the three error terms, the last provides the tightest limit on the size of $L.$ Therefore, any attempt to improve the result in Theorem \ref{Main Result} should begin with refining the proof of Lemma \ref{Log error term} whence this error term comes. However, the primary area for further work is to extend the techniques developed in this paper to the problem of proving an analog for \cite[Theorem 2.2]{SS25} determining how small a murmuration sum can be used to replace the sum over all $n \in \mathbb{N}.$

\subsection*{Acknowledgments}

This paper's first draft was completed in 2026 as a junior paper in partial fulfillment of the requirements
for the Bachelor's degree in mathematics at Princeton University. The author would like to acknowledge his advisor Will Sawin for his guidance, feedback, and innumerable technical suggestions; secondly, Tristen Brisky and Aiden Sharp for useful edits, and suggestions; and finally Luna Dole for helpful conversations regarding Lemma \ref{Improved Trace Bound}. The author assumes sole responsibility for errors in this paper.

\subsection*{AI Declaration} This paper was written without the use of artificial intelligence.

\section{Proofs}

\subsection{Preliminary Results}

We begin by recalling a necessary preliminary result, originally \cite[Lemma 11]{SS25}.

\begin{lemma} \label{Voronoi Summation}
    Let $E$ be an elliptic curve over $\mathbb{Q}$, $N(E)$ the conductor of $E,$ and $\epsilon(E)$ the epsilon factor of $E.$ Let $q$ be a positive integer, $a$ an integer coprime to $q,$ and $W: (0, \infty) \rightarrow \infty) \rightarrow \mathbb{R}$ a smooth compactly supported function. Let $\overline{a N(E)}$ denote the modular inverse of $aN(E)$ modulo $q.$ We have
    $$\frac{\epsilon(E)}{N(E)} \sum_{n = 1}^\infty a_n(E) W\left(  \frac{n}{N(E)} \right) e \left( \frac{an}{q}  \right)$$
    $$= \frac{1}{q} \sum_{n = 1}^{\infty} \frac{a_n(E)}{\sqrt{n}} e \left( \frac{\overline{aN(E)}n}{q}  \right) \int_0^\infty 2\pi\sqrt{u} W(u) J_1\left( 4\pi \frac{\sqrt{un}}{q} \right) du.$$
\end{lemma}

Here, $J_1$ denotes a Bessel function. Note that the formulation of this lemma in \cite{SS25} is erroneous since it contains an extra factor of $\sqrt{\frac{n}{N(E)}}$ in the summation on the left hand side. Lemma \ref{Voronoi Summation} is proved through an application of Voronoi summation, as formulated in \cite[Lemma 2.21]{BFKMMS23} (also erroneously stated, as explained in \cite{SS25}).

\begin{lemma}\label{Test function integral bound}
    Let $W: (0, \infty) \rightarrow \mathbb{R}$ be a smooth, compactly supported function. Let $m \geq 0$, $n > 0,$ and $x > 0$ be  real numbers. For fixed $x > 0$ and sufficiently large varying $L > L_{0}$ (with $L_0$ dependent on $W$ and $n$)
    $$\int_0^\infty W\left( L(u - x) \right) 2\pi \sqrt{u} J_1(4\pi \sqrt{un}) du = O_{m, W} \left( L^{m - 1} n^{-\frac{m}{2} - \frac{1}{4}}\right).$$
\end{lemma}

This lemma is a modification of \cite[Lemma 12]{SS25} to account for the extra parameter $L,$ though the factor of $n$ is also improved by a power of $\frac{1}{4}$. Additionally, here $m$ is permitted to be a real number rather than just an integer. Note that these bounds do not account for dependence on $x,$ which is considered fixed.

\begin{proof}
    Let $Y$ be the smooth, compactly supported function such that $Y(u) = W(L(u - x)).$ We rely on the proof of \cite[Lemma 12]{SS25}, which uses the identity $t J_1(t) = \left( \frac{d}{tdt}  \right)^k (t^{k + 1} J_{k + 1} (t))$ \cite[8.472]{GR14} and iterated integration by parts to obtain the relation
    $$\int_0^\infty Y(u) 2\pi \sqrt{u} J_1 (4\pi \sqrt{un}) du = \frac{(-1)^k}{2\sqrt{n}} \int_0^\infty \frac{d^k Y(u)}{du^k} \cdot \frac{(4\pi \sqrt{un})^{k + 1}}{(8\pi n)^k} J_{k + 1}(4\pi \sqrt{un}) du.$$
    Because $|J_{k + 1}(z)| \leq O \left( \min\left(1, \frac{1}{\sqrt{|z|}} \right) \right),$ for any integer $n,$ we have $J_{k + 1}(4\pi \sqrt{un}) \leq O(1)$ for $0 < u < \frac{1}{n}$ and $J_{k + 1}(4\pi \sqrt{un}) \leq O\left(u^{-\frac{1}{4}} n^{-\frac{1}{4}} \right)$ for $u \geq \frac{1}{n}$. Thus, our expression is
    $$\ll_{k, W} \frac{1}{\sqrt{n}} \int_0^{\frac{1}{n}} \left| \frac{d^k Y(u)}{du^k}   \right| \cdot  \frac{u^\frac{k + 1}{2} \cdot n^{\frac{k + 1}{2}}}{n^k} du + \frac{1}{\sqrt{n}} \int_{\frac{1}{n}}^\infty \left| \frac{d^k Y(u)}{du^k}   \right| \cdot  \frac{u^{\frac{k + 1}{2} - \frac{1}{4}} \cdot n^{\frac{k + 1}{2} - \frac{1}{4}}}{n^k} du.$$ 
    We have that
    $$\frac{d^k Y(u)}{du^k} = \frac{d^k W(L(u - x))}{du^k} = L^k \cdot W^{(k)}(L(u - x)),$$
    so performing the substitution of variables $v = L(u - x)$ gives
    $$L^{k - 1} n^{-\frac{k}{2} }\int_{-Lx}^{\frac{L}{n} - Lx} \left| W^{(k)}(v) \right| \left(\frac{v}{L} + x \right)^\frac{k + 1}{2} dv + L^{k - 1} n^{-\frac{k}{2} - \frac{1}{4}}\int_{\frac{L}{n} -Lx}^{\infty} \left| W^{(k)}(v) \right| \left(\frac{v}{L} + x \right)^{\frac{k + 1}{2} - \frac{1}{4}} dv.$$
    We can factor out $L^{- \frac{k + 1}{2}}$ and $L^{- \frac{k + 1}{2} + \frac{1}{4}}$ from the respective integrals; furthermore, using the fact that $W$ (and thus $W^{(k)}$) has some compact support $[A, B]$, we upper-bound by
    $$L^{\frac{k - 3}{2}} n^{-\frac{k}{2} }\int_{-Lx}^{\frac{L}{n} - Lx} \left| W^{(k)}(v) \right| \left(v + Lx \right)^\frac{k + 1}{2} dv + L^{\frac{k - 3}{2} + \frac{1}{4}} n^{-\frac{k}{2} - \frac{1}{4}}\int_{A}^{B} \left| W^{(k)}(v) \right| \left(v + Lx \right)^{\frac{k + 1}{2} - \frac{1}{4}} dv.$$
    To bound the first integral, assume that $x \leq \frac{1}{n}.$ Otherwise, we have $\frac{L}{n} - Lx < A$ for large $L,$ so the integral evaluates to zero. With our assumption, we get
    \begin{align*}
        L^{\frac{k - 3}{2}} n^{-\frac{k}{2} }\int_{-Lx}^{\frac{L}{n} - Lx} \left| W^{(k)}(v) \right| \left(v + Lx \right)^\frac{k + 1}{2} dv &\leq L^{\frac{k - 3}{2}} n^{-\frac{k}{2} } \sup_v \left| W^{(k)}(v) \right| \int_{A}^{B} \left(v + Lx \right)^\frac{k + 1}{2} dv\\
        &\ll_{k, W} L^{\frac{k - 3}{2}} n^{-\frac{k}{2} } \cdot (B - A) (Lx)^{\frac{k + 1}{2}} \\
        &\ll_{k, W} L^{k - 1} n^{-\frac{k}{2} - \frac{1}{4}} \qquad \left(\text{since $x \leq \frac{1}{n}$}\right)
    \end{align*}
    as desired. For the second integral, we bound
    \begin{align*}
        L^{\frac{k - 3}{2} + \frac{1}{4}} n^{-\frac{k}{2} - \frac{1}{4}}\int_{A}^{B} \left| W^{(k)}(v) \right| \left(v + Lx \right)^{\frac{k + 1}{2} - \frac{1}{4}} dv &\ll_{k, W} L^{\frac{k - 3}{2} + \frac{1}{4}} n^{-\frac{k}{2} - \frac{1}{4}}\int_{A}^{B}  \left(v + Lx \right)^{\frac{k + 1}{2} - \frac{1}{4}} dv\\
        &\ll_{k, W} L^{\frac{k - 3}{2} + \frac{1}{4}} n^{-\frac{k}{2} - \frac{1}{4}} \cdot (B - A) (Lx)^{\frac{k + 1}{2} - \frac{1}{4}}.\\
        &\ll_{k, W} L^{k - 1} n^{-\frac{k}{2} - \frac{1}{4}}.
    \end{align*}
    The above holds in particular for $k = \lfloor m\rfloor$ and $k = \lceil m \rceil.$ Since we can choose some positive $a, b$ with $a + b = 1$ satisfying $m = a\lfloor m \rfloor +  b\lceil m \rceil$, taking a weighted geometric mean of the bounds $O_{k, W} \left( L^{k - 1} n^{-\frac{k}{2} - \frac{1}{4}} \right)$ for $k = \lfloor m \rfloor$ and $k = \lceil m \rceil$ establishes the bound for $m \in \mathbb{R}_{\geq 0}$.
\end{proof}

\begin{defn}\label{def q}
    To each integer $n$, we associate another integer
$$q(n) := 2^{d_2(n)} \cdot 3^{d_3(n)} \cdot \prod_{\substack{p \mid n, p \neq 2, 3}} p,$$
where $d_p(n)$ is $7$ if $p \mid n$ and $0$ otherwise.
\end{defn}

\begin{lemma} \label{dependence of coefficients on $q$}
    Let $q(n)$ be as defined as above. Let $a, b, a', b' \in \mathbb{Z}$. Suppose that there are no primes $p$ such that both $p^4 \mid a$ and $p^6 \mid b$ and no primes with $p^4 \mid a'$ and $p^6 \mid b'.$ Then
    $$a \equiv a', b \equiv b' \pmod {q(n)} \implies a_n(E_{a, b}) = a_n(E_{a', b'}).$$
\end{lemma}

\begin{proof}
    By the multiplicativity of $a_n$ and $q(n),$ it suffices to prove the result for $n = p^v.$ If $p \neq 2, 3$, the fact that either $p^4 \nmid a$ or $p^6 \nmid b$ holds guarantees that $E_{a, b}$ is minimal and can be reduced. Since $a_n(E_{a, b})$ depends on the number of solutions and reduction type of the reduced curve over $\mathbb{F}_p$ and changing $a$ and $b$ by multiples of $p$ does not change the reduced curve, the result immediately follows. 
    
    If $p = 3,$ then the curve $y^2 = x^3 + ax + b$ is not necessarily a minimal Weierstrass equation. \cite[Lemma 17]{SS25} gives that for any elliptic curve $E = E(\mathbb{F}_3),$ the curve $E_{a, b}$ has reduction $E$ in $\mathbb{F}_3$ if and only if there exists $\lambda, r \in \mathbb{Q}_3,$ $\lambda$ a power of $3$ such that
    $$y^2 = (x + r)^3 + \lambda^4 a(x + r) + \lambda^6 b$$ is a polynomial over $\mathbb{Z}_3$ and reduces to $E$ mod $3.$ Following the proof of \cite[Lemma 16]{SS25}, we start with the case where $E$ is supersingular, in which case the $x^2$ coefficient of $(x + r)^3 + \lambda^4 a(x + r) + \lambda^6 b$ must be $0 \pmod 3$. Then $v_3(r) \geq 0,$ which implies that $v_3(\lambda^4 a) \geq 0$ and $v_3(\lambda^6 b) \geq 0$. But at least one of $v_3(\lambda^4 a)$ or $v_3(\lambda^6 b)$ must be zero, because the reduction $y^2 = (x - r)^3$ would otherwise not be an elliptic curve. Since we cannot have $3^4 \mid a$ and $3^6 \mid b,$ this forces $\lambda = 1.$ Hence, the curve $y^2 = x^3 + ax + b$ is already minimal, and changing $a$ or $b$ by a multiple of $3$ does not change the reduction.

    If $E$ is not supersingular, then the $x^2$ coefficient is not $0$ in $\mathbb{F}_p,$ which forces $v_3(r) = -1.$ Since the coefficients must satisfy
    $$3r^2 + \lambda^4 a \in \mathbb{Z}_3 \quad \text{and} \quad r^3 + \lambda^4 a r + \lambda^6 b \in \mathbb{Z}_3,$$
    we must have $v_3(\lambda^4 a) = -1$, which forces $v_3(\lambda^6 b) \geq -3$ as well. Since we cannot have both $3^4 \mid a$ and $3^6 \mid b,$ we must have $\lambda = 3^{-1}.$ Then, if we add $3^7 a'$ to $a$ and $3^7 b'$ to $b,$ the new curve must be isomorphic to
    $$y^2 = (x + r)^3 + \lambda^4 a(x + r) + 3^3 a'(x + r) + \lambda^6 b + 3b',$$
    which reduces down to the same curve. Therefore, $E_{a, b}$ and $E_{a + 3^7 a', b + b^7 b'}$ are isomorphic.

    Finally, we consider the case of $p = 2.$ \cite[Lemma 22]{SS25} gives that for $E \in \overline{\mathcal{M}}_{1, 1}(\mathbb{F}_2),$ a curve $E_{a, b}$ has reduction mod $2$ isomorphic to $E$ if and only if there exists $r, s, t \in \mathbb{Q}_2$ and $\lambda$ an integer power of $2$ such that
    $$(y + sx + t)^2 = (x + r)^3 + \lambda^4 a(x + r) + \lambda^6 b$$
    has coefficients in $\mathbb{Z}_2$ and reduces mod $2$ to $E.$ Through the proof of the $p = 2$ case of \cite[Lemma 16]{SS25}, we must have either one of $v_2(s) = -1$ and $v_2(t) = -1$, or the $xy$ and $y$ coefficients of the reduced curve would be zero and the curve would be singular. We also have that $v_2(r) \geq -2$ and
    $$-4 \leq v_2(\lambda^4 a) < 0, \qquad -6 \leq v_2(\lambda^4 b) < 0.$$
    Since we cannot have both $2^4 \mid a$ and $2^6 \mid b,$ we must have $\lambda = 2^{-1}$. Thus, similar to the $p = 3$ case, we can be sure that adding multiples of $2^7$ to $a$ and $b$ will not change the reduction of the curve $(y + sx + t)^2 = (x + r)^3 + \lambda^4 a(x + r) + \lambda^6 b.$
\end{proof}
    
\begin{rem}\label{p adic coefficients}
    The proof of Lemma \ref{dependence of coefficients on $q$} demonstrates that the function $(a, b) \rightarrow a_{p^v}(E_{a, b})$ is continuous on $\mathbb{Z}_p^2 \setminus (p^4 \mathbb{Z}_p \times p^6 \mathbb{Z}_p)$ and that for $(a, b) \in \mathbb{Z}_p^2 \setminus (p^4 \mathbb{Z}_p \times p^6 \mathbb{Z}_p),$ we have
    $$E_{p^v(n)}(E_{a, b}) = E_{p^{v_p(n)}}(E_{\overline{a}, \overline{b}}),$$
    where $\overline{a}, \overline{b} \in \mathbb{F}_{p^{v_p(q(n))}}$ are the reductions of $a$ and $b.$ Meanwhile, for $(a, b) \in (p^4 \mathbb{Z}_p \times p^6 \mathbb{Z}_p),$ we simply have $E_{p^{v_p(n)}}(E_{\overline{a}, \overline{b}}) = 0.$
\end{rem}

\begin{lemma}\label{Trace Function Bound}
    Let $q = q(n).$ Then for any $y, z \in \mathbb{Z}/q,$
    $$ \left|\sum_{a, b \in \mathbb{Z}/q}  a_n(E_{a, b}) \cdot e\left( \frac{a y + b z}{q}  \right) \right| \leq O_\epsilon\left( q^{\frac{3}{2} + \epsilon} \cdot n^{\frac{1}{2} + \epsilon}\right).$$
\end{lemma}

\begin{proof}
    We begin by factoring the sum as
    \begin{align*}
        \left|\sum_{a, b \in \mathbb{Z}/q}  a_n(E_{a, b}) \cdot e\left( \frac{a y + b z}{q}  \right) \right| &= \left| \sum_{a, b \in \mathbb{Z}/q}  \prod_{p \mid q}  \left( a_{p^{v_p(n)}}(E_{a, b}) \cdot e\left( \frac{a y_p + b z_p}{p^{v_p(q)}}  \right) \right) \right|\\
        &= \prod_{p \mid q} \left| \sum_{a, b \in \mathbb{Z}/p^{v_p(q)}}    a_{p^{v_p(n)}}(E_{a, b}) \cdot e\left( \frac{a y_p + b z_p}{p^{v_p(q)}}  \right)  \right|,
    \end{align*}
    by the Chinese remainder theorem; here, we define $y_p = y \cdot (q/p^{v_p(q)})^{-1} \in \mathbb{Z}/p^{v_p(q)}$ and $z_p$ is defined likewise. Recall that $a_{p^{v_p(n)}}(E_{a, b})$ is well defined by Lemma \ref{dependence of coefficients on $q$}. Recall the definitions and bounds
    \begin{align}
        \sigma_0(n) &:= \#\{ k \in \mathbb{N} : k \mid n\} \ll \sigma_0(n) \leq 2^{1.53794 \, \frac{\log n}{\log \log n}} & \text{by \cite[Th\'eor\`eme 1]{NR83}}\label{eqn:sigma-bound},\\
        \omega(q) &:= \#\{ p \text{ prime} : p \mid n\} \leq 1.38407 \, \frac{\log n}{\log \log n} & \text{by \cite[Th\'eor\`eme 11]{Rob83}}\label{eqn:omega-bound}.
    \end{align}

    It suffices to show there exists a constant $C$ independent of $n$ and $q$ such that
    \begin{equation*}
        \sum_{a, b \in \mathbb{Z}/p^{v_p(q)}} a_{p^{v_p(n)}}(E_{a, b}) \cdot e\left( \frac{a y + b z}{p^{v_p(q)}}  \right) < C \cdot p^{\frac{3}{2}} \cdot \sigma_0\left( p^{v_p(n)} \right) \sqrt{p^{v_p(n)}} \qquad \text{for any $y, z \in \mathbb{Z}/p^{v_p(q)}$}
    \end{equation*}
    since multiplying over all $p \mid q$ gives
    $$\left|\sum_{a, b \in \mathbb{Z}/q}  a_n(E_{a, b}) \cdot e\left( \frac{a y + b z}{q}  \right) \right| \ll_\epsilon C^{\omega(q)} \cdot q^{\frac{3}{2}} \cdot \sigma_0(n) n^{\frac{1}{2}} \ll_\epsilon q^{\frac{3}{2} + \epsilon} \cdot n^{\frac{3}{2} + \epsilon},$$
    with the last inequality by \eqref{eqn:sigma-bound} and \eqref{eqn:omega-bound}. For $p \in \{2, 3\},$ the Hasse-Weil bound gives
    $$\sum_{a, b \in \mathbb{Z}/p^{v_p(q)}} a_{p^{v_p(n)}}(E_{a, b}) \cdot e\left( \frac{a y + b z}{p^{v_p(q)}}  \right) \leq p^{14} \cdot \sigma_0 \left(p^{v_p(n)} \right)  \sqrt{p^{v_p(n)}},$$
    so it suffices to take $C$ at least some high enough power of $3$. For $p \neq 2, 3$, we first consider the case where at least one of the following conditions is true: $(y, z) \neq (0, 0)$ or $v$ is odd. For each pair $(a, b) \in (\mathbb{Z}/p)^\times \times (\mathbb{Z}/p)^\times,$ there exists a unique pair $(c, t) \in (\mathbb{Z}/p)^\times \times (\mathbb{Z}/p)^\times$ such that $(a, b) = (ct^2, ct^3).$ This follows since the map $(c, t) \rightarrow (ct^2, ct^3)$ has inverse $(a, b) \rightarrow \left( \frac{a^3}{b^2}, \frac{b}{a} \right).$ This allows us to rewrite the $p$-factor of the desired product as
    \begin{align*}
        a_{p^v}(E_{0, 0}) &+ \sum_{a \in (\mathbb{Z}/p)^\times} a_{p^v}(E_{a, 0}) e \left( \frac{ay}{p} \right) + \sum_{b \in (\mathbb{Z}/p)^\times} a_{p^v}(E_{0, b}) e \left( \frac{bz}{p} \right) \\
        &+ \sum_{c, t \in (\mathbb{Z}/p)^\times}a_{p^v}(E_{ct^2, ct^3}) e \left( \frac{ct^2 y + ct^3 z}{p} \right).
    \end{align*}
    Hasse-Weil bounds the first three summations by 
    $$\left| a_{p^v}(E_{0, 0}) + \sum_{a \in (\mathbb{Z}/p)^\times} a_{p^v}(E_{a, 0}) e \left( \frac{ay}{p} \right) + \sum_{b \in (\mathbb{Z}/p)^\times} a_{p^v}(E_{0, b}) e \left( \frac{bz}{p} \right) \right| \leq 4 p^{1 + \frac{v}{2}} \sigma\left(p^{\frac{v}{2}} \right) \ll 4 p^{\frac{3}{2} + \frac{v}{2}}.$$
    For the fourth summation, we recall the standard identity (see, e.g., \cite[Exercise 4.10]{washington2008elliptic})
    \begin{equation}
        a_{p^v}(E_{t^2 a, t^3 b}) = a_{p^v}(E_{a, b}) \lg{t}{p}^v \qquad \text{for $t \in (\mathbb{Z}/p)^\times$} \label{eqn:quadratic-twist-lemma}
    \end{equation}
    to obtain
    
    \begin{align*}
        \sum_{c, t \in (\mathbb{Z}/p)^\times}a_{p^v}(E_{ct^2, ct^3}) e \left( \frac{ct^2 y + ct^3 z}{p} \right) &= \sum_{c \in (\mathbb{Z}/p)^\times}a_{p^v}(E_{c, c}) \sum_{(t \in \mathbb{Z}/p)^\times} e \left( \frac{ct^2 y + ct^3 z}{p} \right) \lg{t}{p}^v.
    \end{align*}
    We employ the the bound
    $$\left| \sum_{x} \psi\left[ F(x) \right] \chi\left[ R_\mathfrak{d}(x) \right] \right| \leq (k + d - 1) \sqrt{q},$$
    included as the last equation of \cite{W48}. Here, $\psi$ is a nontrivial character of $(\mathbb{F}_q, +)$ and $\chi$ is a nontrivial character of $\mathbb{F}_q^\times,$ while $F$ is a polynomial of degree $k$ with constant term $0$ and $R_\mathfrak{d}$ is a degree $d$ polynomial. The bound requires either than $R_\mathfrak{d}$ not to be a constant times the $\text{ord}(\chi)$-th power of a polynomial, or that $F$ is not of the form $g^p - g$ for some polynomial $g.$ We let
    $$\psi = e\left( \frac{\cdot}{p} \right), \quad \chi = \lg{\cdot}{p}, \quad F(t) = ct^2y + ct^3z, \quad R_\mathfrak{d}(t) = t^v.$$
    If $v$ is odd, then $R_\mathfrak{d}(t) = t^v$ is not a constant times a $2$nd power, while if $(y, z) \neq 0,$ then $F(t) = ct^2y + ct^3 y$ is never of the form $g^p - g$ for a polynomial $g$. Therefore, our assumptions on $x, y,$ and $v$ guarantee that the non-degeneracy conditions hold. We apply the bound and obtain
    $$\left| \sum_{(t \in \mathbb{Z}/p)^\times} e \left( \frac{ct^2 y + ct^3 z}{p} \right) \lg{t}{p}^v \right| \leq 4 \sqrt{p}.$$
    Then a simple application of the Hasse-Weil bound suffices to conclude
    $$\left| \sum_{c \in (\mathbb{Z}/p)^\times}a_{p^v}(E_{c, c}) \sum_{(t \in \mathbb{Z}/p)^\times} e \left( \frac{ct^2 y + ct^3 z}{p} \right) \lg{t}{p}^v \right| \leq 4\sqrt{p} \sum_{c \in (\mathbb{Z}/p)^\times} \left|a_{p^v}(E_{c, c}) \right| \leq_\epsilon 8p^{\frac{3}{2}} \cdot p^{\frac{v}{2} + \epsilon}.$$
    since $v \geq 2.$ Combined with the error of $4 p^{\frac{3}{2}} \cdot p^{\frac{v}{2}},$ our total term for $p$ is less than $12 p^{\frac{3}{2}} \cdot p^{\frac{v}{2} + \epsilon},$

    In the case where $y = z = 0$ and $v$ is even, then \eqref{eqn:quadratic-twist-lemma} gives that $a_{p^v}(E_{a, b}) = a_{p^v}(E_{at^2, bt^3})$ for any $t \in (\mathbb{Z}/p)^\times.$ The number of $t$ such that $(a, b) = (t^4 a, t^6 b)$ is equal to $|\text{Aut}(E)|.$ Therefore, each isomorphism class of elliptic curves $E$ corresponds to $\frac{p - 1}{|\text{Aut}(E)|}$ pairs $(a, b)$, so
    $$\sum_{a, b \in \mathbb{Z}/p} a_{p_v}(E_{a, b}) = a_{p_v}(E_{0, 0}) + \sum_{E} a_{p^v}(E)\frac{p - 1}{|\text{Aut(E)}|} = - (p - 1) \sum_{\substack{f \in S_0^{v + 2}(\text{SL}_2(\mathbb{Z})) \\ \text{eigenform} \\ a_1(f) = 1}} a_p(f),$$
    by \cite[Lemma 25]{SS25} (noting that $a_{p_v}(E_{0, 0}) = 0$). But we have that $|a_p(f)| \leq 2p^{\frac{v + 1}{2}}$ by the Ramanujan-Petersson conjecture for holomorphic cusp forms, as proved by Deligne in \cite{D68} and \cite{D74}. For even $v$, the number of $f$ in the sum is bounded by 
    $$\dim S_0^{v + 2} \leq \dim M_{v + 2} - 1 \leq \lfloor (v + 2)/12 \rfloor + 1 - 1 \leq \frac{v + 2}{12},$$
    where the bound $\dim M_{2k} \leq \lfloor k/6 \rfloor$ is standard and found, for instance, in \cite[Corollary VII.3.2.1]{S73}. Thus, we conclude that the total sum is upper bounded by $\frac{v + 2}{6} \cdot 2p^{\frac{v + 3}{2}}.$ We note that
    $$\prod_{p \neq 2, 3} \frac{v_p(n) + 2}{6} < \prod_{p \neq 2, 3} (v_p(n) + 1) < \sigma_0(n) \leq O_\epsilon(n^\epsilon),$$
    where the final inequality is due to \eqref{eqn:sigma-bound}. This guarantees that the factors $\frac{v + 2}{6} \cdot 2p^{\frac{v + 3}{2}}$ (along with the factors for $p = 2$ and $3$) multiply to a product at most $O_\epsilon\left( q^{\frac{3}{2}} n^{\frac{1}{2} + \epsilon} \right)$.
\end{proof}

In fact, the bound in Lemma \ref{Trace Function Bound} can be improved by using a different technique for primes $p$ with $v_p(n) = 1.$

\begin{lemma} \label{Improved Trace Bound}
    Let $q = q(n).$ Then
    $$\left|\sum_{a, b \in \mathbb{Z}/q}  a_n(E_{a, b}) \cdot e\left( \frac{a y + b z}{q}  \right) \right| \leq O_\epsilon \left(  q^{\frac{3}{2} + \epsilon}\prod_{p : p^2 \mid n} p^{\frac{v_p(n)}{2} + \epsilon} \right).$$
\end{lemma}

\begin{rem}
    Lemma \ref{Improved Trace Bound} is sufficient to prove that the sum in question is always $\ll_\epsilon n^{\frac{3}{2} + \epsilon}.$ However, this particular fact is not useful for our purposes.
\end{rem}

\begin{proof}
    It suffices to amend the proof of Lemma \ref{Trace Function Bound} by proving that for primes $p$ with $v_p(n) = 1$ and a large enough absolute constant $C,$
    $$\left| \sum_{a, b \in \mathbb{Z}/p^{v_p(q)}} a_{p^{v_p(n)}}(E_{a, b}) e\left( \frac{a y + b z}{p^{v_p(q)}}  \right) \right| \leq C p^{\frac{3v_p(q)}{2}}.$$
    The cases of $p \in \{2, 3\}$ once again follow from an application of the Hasse-Weil bound. For $p \notin \{2, 3\},$ first suppose that $z \neq 0$. Then we can write
    \begin{align}
        &\quad \sum_{a, b \in \mathbb{Z}/p} a_{p}(E_{a, b}) e\left( \frac{a y + b z}{p}  \right) \nonumber\\
        &= -\sum_{x \in \mathbb{Z}/p} \sum_{a \in \mathbb{Z}/p} \sum_{b \in \mathbb{Z}/p} \lg{x^3 + ax + b}{p} e\left( \frac{a y + b z}{p}  \right) \nonumber\\
        &= -\sum_{x \in \mathbb{Z}/p} \sum_{a \in \mathbb{Z}/p} \sum_{b' \in \mathbb{Z}/p} \lg{b'}{p} e\left( \frac{a y + b' z - (x^3 - ax)z}{p}  \right) \label{eqn:cancelling-Legendre-sum}\\
        & \qquad \text{after making the substitution $b' := x^3 + ax + b$; $\lg{\cdot }{p}$ denotes the Legendre symbol} \nonumber\\
        &= -\sqrt{ \lg{-1}{p} p} \cdot \sum_{x \in \mathbb{Z}/p} \sum_{a \in \mathbb{Z}/p} e\left( \frac{a y - (x^3 - ax)z}{p}  \right) \nonumber\\
        &\qquad \text{by the Gauss sum identity $\sum_{b' \in \mathbb{Z}/p} \lg{b'}{p} e \left( \frac{b'}{p}  \right) = \sqrt{ \lg{-1}{p} p}$}\nonumber\\
        &= -p \cdot \sqrt{ \lg{-1}{p} p} \cdot e\left( -\frac{y^3/z^2}{p}\right), \qquad \qquad \left( \text{since $\sum_{a \in \mathbb{Z}/p} e\left( \frac{a y + axz}{p}  \right) = \begin{cases}
        p & y = -xz\\
        0 & \text{otherwise}
    \end{cases}$} \right) \nonumber
    \end{align}
    which is $p^{\frac{3}{2}}$ in absolute value. If $z = 0$, then the sum \eqref{eqn:cancelling-Legendre-sum} cancels to zero.
\end{proof}

\subsection{Averages of \texorpdfstring{$L$}{L}-function Coefficients}

For the remainder of this paper, we will work in the following setting. Let $X$ be a fixed positive real number, let $n$ be any natural number, and let $q := q(n)$ as in Definition \ref{def q}. In order to prove Theorem \ref{Main Result}, we will need to prove the following equality.

\begin{prop} \label{sums of coefficients}
    Let $q = q(n).$ Then, as $X$ and $n$ grow,
    \begin{align*}
        \sum_{\substack{A, B \in \mathbb{Z} \\ \max(4|A|^3, 27|B|^2) \leq X \\ \text{$p^4 \nmid A$ or $p^6 \nmid B$ for all $p$}\\ 4A^3 + 27 B^2 \neq 0}} \! \! \! \! \! \! \! \! \! \! \! \! \! \! &a_n(E_{A, B}) = \frac{1}{q^2} \cdot 2(X/4)^{\frac{1}{3}} \cdot 2(X/27)^{\frac{1}{2}} \cdot \prod_{p \nmid q} \left( 1 - \frac{1}{p^{10}} \right) \cdot \sum_{a, b \in \mathbb{Z}/q} a_n(E_{a, b}) \\
        &+ O_\epsilon \left( X^{\frac{1}{6}} n^{\frac{1}{2} + \epsilon} \right)+ O_\epsilon\left(X^{\frac{1}{2}} q^{\frac{1}{2} + \epsilon}\prod_{p : p^2 \mid n} p^{\frac{v_p(n)}{2} + \epsilon} \right) + O_\epsilon \left( X^{\frac{1}{12}} q^{\frac{3}{2} + \epsilon}\prod_{p : p^2 \mid n} p^{\frac{v_p(n)}{2} + \epsilon} \right).
    \end{align*}
    
\end{prop}
Proposition \ref{sums of coefficients} is important because it will allow us to average the coefficient $a_n$ over all elliptic curves of naive height less than or equal to $X$ and thus prove Lemma \ref{Average of coefficients}. The first step in proving Proposition \ref{sums of coefficients} is the following reformulation of the left hand side.

\begin{lemma}\label{Polya Vinogradov}
    Let $q = q(n).$ Then, as $X$ and $n$ grow,
    \begin{align*}
       \sum_{\substack{A, B \in \mathbb{Z} \\ \max(4|A|^3, 27|B|^2) \leq X \\ \text{$p^4 \nmid A$ or $p^6 \nmid B$ for all $p$}\\ 4A^3 + 27 B^2 \neq 0}} &a_n(E_{A, B}) = \frac{1}{q^2} \sum_{y, z \in \mathbb{Z}/q} \left( \sum_{a, b \in \mathbb{Z}/q}  a_n(E_{a, b}) \cdot e\left( \frac{a y + b z}{q}  \right)\right)\\
       &\cdot  \sum_{\substack{A, B \in \mathbb{Z} \\ \max(4|A|^3, 27|B|^2) \leq X \\ \text{$p^4 \nmid A$ or $p^6 \nmid B$ for all $p$}}}  e \left( \frac{-A y - Bz}{q}  \right) + O_\epsilon \left( X^{\frac{1}{6}} n^{\frac{1}{2} + \epsilon} \right).
    \end{align*}
\end{lemma}
\begin{proof}
    The number of $A, B$ with $\max(4|A|^3, 27|B|^2) \leq X$ such that $4A^3 + 27B^2 = 0$ is at most $O(X^{\frac{1}{6}}).$ Therefore, the Hasse-Weil bound suffices to show that we get the error term by dropping the condition $4A^3 + 27B^2 \neq 0.$
    We distribute the product of sums on the right hand side and rewrite the main term as a product over $a, b, A,$ and $B$:
    \begin{align*}
        & \quad \frac{1}{q^2} \sum_{y, z \in \mathbb{Z}/q} \sum_{\substack{a, b \in \mathbb{Z}/q\\ A, B \in \mathbb{Z} \\ \max(4|A|^3, 27|B|^2) \leq X \\ \text{$p^4 \nmid A$ or $p^6 \nmid B$ for all $p$}}}  a_n(E_{a, b}) \cdot e\left( \frac{ay + bz}{q}   \right) \cdot e\left( \frac{-Ay - Bz}{q}  \right)\\
        &= \frac{1}{q^2} \sum_{\substack{a, b \in \mathbb{Z}/q\\ A, B \in \mathbb{Z} \\ \max(4|A|^3, 27|B|^2) \leq X \\ \text{$p^4 \nmid A$ or $p^6 \nmid B$ for all $p$}}}  a_n(E_{a, b}) \sum_{y, z \in \mathbb{Z}/q} e\left( \frac{ay + bz - Ay - Bz}{q}   \right)\\
        &= \sum_{y, z \in \mathbb{Z}/q} \sum_{\substack{a, b \in \mathbb{Z}/q\\ A, B \in \mathbb{Z} \\ \max(4|A|^3, 27|B|^2) \leq X \\ \text{$p^4 \nmid A$ or $p^6 \nmid b$ for all $p$}}}  a_n(E_{a, b})\\
        &= \sum_{\substack{A, B \in \mathbb{Z} \\ \max(4|A|^3, 27|B|^2) \leq X \\ \text{$p^4 \nmid A$ or $p^6 \nmid b$ for all $p$}}}  a_n(E_{a, b}),
    \end{align*}
    where the penultimate step is due to the identity
    \begin{align*}
        \sum_{y, z \in \mathbb{Z}/q} e\left( \frac{ay + bz - Ay - Bz}{q}   \right) &= \left( \sum_{y \in \mathbb{Z}/q} e\left( \frac{(a - A)y}{q}  \right)  \right) \left( \sum_{z \in \mathbb{Z}/q} e\left( \frac{(b - B)z}{q}  \right) \right)\\
        &= \begin{cases}
            q^2 & a - A \equiv b - B \equiv 0 \pmod q\\
            0 & \text{otherwise}.
        \end{cases} \qedhere
    \end{align*}
\end{proof}

The next lemma will contribute the main term in Proposition \ref{sums of coefficients}.

\begin{lemma} \label{Main Term}
    Let $q = q(n).$ Then, as $X$ and $n$ grow,
    $$\frac{1}{q^2} \sum_{y, z \in \mathbb{Z}/q} \left( \sum_{a, b \in \mathbb{Z}/q}  a_n(E_{a, b}) \cdot e\left( \frac{a y + b z}{q}  \right)\right)\left( \sum_{\substack{m \leq X^{\frac{1}{12}}\\ q \mid m^4 y\\ q \mid m^6 z}} \mu(m) \sum_{\substack{A, B \in \mathbb{Z} \\ \max(4|A|^3, 27|B|^2) \leq X \\ m^4 \mid A, m^6 \mid B}}  e \left( \frac{-A y - Bz}{q}  \right) \right)$$
    $$= \frac{1}{q^2} \cdot 2(X/4)^{\frac{1}{3}} \cdot 2(X/27)^{\frac{1}{2}} \cdot \prod_{p \nmid q} \left( 1 - \frac{1}{p^{10}} \right) \cdot \sum_{a, b \in \mathbb{Z}/q} a_n(E_{a, b}) + O_\epsilon\left( X^{\frac{1}{2}} q^{ - \frac{1}{2} + \epsilon}\prod_{p : p^2 \mid n} p^{\frac{v_p(n)}{2} + \epsilon} \right).$$
\end{lemma}

\begin{proof}
    First, we note that for $m$ satisfying $q \mid (m^4 y, m^6 z):$
    \begin{align*}
        \sum_{\substack{A, B \in \mathbb{Z} \\ \max(4|A|^3, 27|B|^2) \leq X \\ m^4 \mid A, m^6 \mid B}}  e \left( \frac{-A y - Bz}{q}  \right)  &= \sum_{\substack{A, B \in \mathbb{Z} \\ \max(4|A|^3, 27|B|^2) \leq X \\ m^4 \mid A, m^6 \mid B}}  1\\ 
        &= \left( \frac{2(X/4)^{\frac{1}{3}}}{m^4} + O(1)  \right) \left( \frac{2(X/27)^{\frac{1}{2}} }{m^6} + O(1)  \right)\\
        &= \frac{2(X/4)^{\frac{1}{3}} \cdot 2(X/27)^{\frac{1}{2}}}{m^{10}} + O\left( \frac{2(X/27)^{\frac{1}{2}} }{m^6} +  \frac{2(X/4)^{\frac{1}{3}}}{m^4}  \right).
    \end{align*}
    Postponing the error terms until later, our main term is
    $$2(X/4)^{\frac{1}{3}} \cdot 2(X/27)^{\frac{1}{2}} \cdot \frac{1}{q^2} \sum_{y, z \in \mathbb{Z}/q} \left( \sum_{a, b \in \mathbb{Z}/q}  a_n(E_{a, b}) \cdot e\left( \frac{a y + b z}{q}  \right)\right)\left( \sum_{\substack{m \leq X^{\frac{1}{12}}\\ q \mid m^4 y\\ q \mid m^6 z}} \frac{\mu(m)}{m^{10}} \right).$$
    Because the condition $q \mid \gcd(m^4y, m^6 z)$ depends only on the part of $m$ that shares prime factors with $q,$ we can rewrite
    $$\sum_{\substack{m \leq X^{\frac{1}{12}}\\ q \mid m^4 y\\ q \mid m^6 z}} \frac{\mu(m)}{m^{10}} = \sum_{\substack{k \mid q\\ q|\gcd (k^4 y, k^6 z)}} \sum_{\substack{m' \leq X^{\frac{1}{12}}/k\\ \gcd(m',q) = 1}} \frac{\mu(m'  k)}{(m' k)^{10}} = \sum_{\substack{k \mid q\\ q|\gcd (k^4 y, k^6 z)}} \frac{\mu(k)}{k^{10}} \sum_{\substack{m' \leq X^{\frac{1}{12}}/k\\ \gcd(m',q) = 1}} \frac{\mu(m')}{(m')^{10}}.$$ 
    We need only sum over $k$ square-free since otherwise $\mu(k)$ just evaluates to zero. Switching the orders of the summations, we find that our main term is equal to $2(X/4)^{\frac{1}{3}} \cdot 2(X/27)^{\frac{1}{2}} \cdot \frac{1}{q^2}$ times the sum
    $$\sum_{\substack{k \mid q\\ \text{$k$ square-free}}} \frac{\mu(k)}{k^{10}} \left(\sum_{\substack{m' \leq X^{\frac{1}{12}}/k\\ \gcd(m',q) = 1}} \frac{\mu(m')}{(m')^{10}} \right) \sum_{a, b \in \mathbb{Z}/q}  a_n(E_{a, b}) \sum_{\substack{y, z \in \mathbb{Z}/q\\ q \mid \gcd (k^4 y, k^6 z)}} e\left( \frac{a y + b z}{q}  \right).$$
    For fixed $k,$ we determine the choices of $y$ and $z$ that satisfy $q \mid \gcd (k^4 y, k^6 z).$ For $p \neq 2, 3$ dividing $q,$ if $p \mid k$, then there are no restrictions on $v_p(y)$ and $v_p(z),$ whereas if $p \nmid k,$ then we require $v_p(y) = v_p(z) = 1$ (recall that these $p$-adic valuations cannot exceed $1$ because $q$ is defined so that $v_p(q) = 1.$ In the case where $2 \mid q$: if $2 \nmid k,$ we require that $v_2(y) = v_2(z) = 7,$ while if $2 \mid k,$ we require $v_2(y) \geq 3$ and $v_2(y) \geq 1$ (since $v_2(q) = 8$). Likewise for $3.$ In total, we are summing precisely over $y \in \mathbb{Z}/q$ divisible by $c_y(k)$ and $z$ divisible by $c_z(k)$ where
    $$c_y(k) := 2^{d_y(2, k)} \cdot 3^{d_y(3, q, k)}\prod_{\substack{p \mid q, p \nmid k \\ p \neq 2, 3}} p , \quad c_z(k) := 2^{d_z(2, k)} \cdot 3^{d_z(3, q, k)}\prod_{\substack{p \mid q, p \nmid k \\ p \neq 2, 3}} p $$
    with 
    $$d_y(p, k) = \begin{cases}
        0 & p \nmid q\\
        3 & m \mid k\\
        7 & p \mid q, p \nmid q
    \end{cases} \quad \text{and} \quad d_z(p, k) = \begin{cases}
        0 & p \nmid q\\
        1 & m \mid k\\
        7 & p \mid q, p \nmid q.
    \end{cases} $$
    Then we can rewrite the sum
    $$\sum_{\substack{y, z \in \mathbb{Z}/q\\ q \mid \gcd (k^4 y, k^6 z)}} e\left( \frac{a y + b z}{q}  \right) = \left( \sum_{c_y(k) \mid y} e\left(  \frac{ay}{q} \right) \right) \left( \sum_{c_z(k) \mid z} e\left(  \frac{bz}{q} \right) \right),$$
    this being a product of two geometric sums. At least one of these sums will evaluate to zero unless $q \mid a c_y(k)$ and $q \mid b c_z(k).$ Thus, we can rewrite the main term as $2(X/4)^{\frac{1}{3}} \cdot 2(X/27)^{\frac{1}{2}} \cdot \frac{1}{q^2}$ times the sum
    $$\sum_{\substack{k \mid q\\ \text{$k$ square-free}}} \frac{\mu(k)}{k^{10}} \left(\sum_{\substack{m' \leq X^{\frac{1}{12}}/k\\ (m',q) = 1}} \frac{\mu(m')}{(m')^{10}} \right) \sum_{\substack{a, b \in \mathbb{Z}/q \\ q| \gcd (a c_y, b c_z)}}  a_n(E_{a, b}) \sum_{\substack{y, z \in \mathbb{Z}/q\\ q \mid \gcd (k^4 y, k^6 z)}} e\left( \frac{a y + b z}{q}  \right).$$
    For fixed $k \neq 1$, there exists a prime $p$ dividing $k$. The condition $q \mid \gcd (k^4 y, k^6 z)$ requires that $p$ divide both $a$ and $b.$ Then, the $L$-function coefficient $a_{p^{v_p(n)}}(E_{a, b}) = 0,$ which means that $a_n(E_{a, b}) = \prod_{p \mid n} a_{p^{v_p(n)}}(E_{a, b}) = 0$ as well. Thus, for all $k \neq 1,$ the entire main term sum evaluates to zero. Our main term is thus given by 
    $$2(X/4)^{\frac{1}{3}} \cdot 2(X/27)^{\frac{1}{2}} \cdot \frac{1}{q^2}\left(\sum_{\substack{m' \leq X^{\frac{1}{12}}\\ \gcd(m',q) = 1}} \frac{\mu(m')}{(m')^{10}} \right) \sum_{\substack{a, b \in \mathbb{Z}/q}}  a_n(E_{a, b}) \sum_{\substack{y, z \in \mathbb{Z}/q\\ q \mid \gcd(y, z)}} e\left( \frac{a y + b z}{q}  \right).$$
    We note that $\sum_{\substack{y, z \in \mathbb{Z}/q\\ q \mid \gcd(y, z)}} e\left( \frac{a y + b z}{q}  \right)$ is clearly equal to $1$, while
    \begin{align*}
        \sum_{\substack{m' \leq X^{\frac{1}{12}}\\ (m',q) = 1}} \frac{\mu(m')}{(m')^{10}} &= \sum_{\substack{m' \mid \gcd(m',q) = 1}} \frac{\mu(m')}{(m')^{10}} - \sum_{\substack{m' > X^{\frac{1}{12}}\\ \gcd(m',q) = 1}} \frac{\mu(m')}{(m')^{10}}\\
        &= \prod_{p \nmid m'} \left(1 - \frac{1}{p^{10}} \right) + O \left( \sum_{\substack{m' > X^{\frac{1}{12}}}} \frac{1}{(m')^{10}} \right) \\
        &= \prod_{p \nmid m'} \left(1 - \frac{1}{p^{10}} \right) + O \left( \frac{1}{X^{\frac{3}{4}}} \right).
    \end{align*}
    Our term thus gives us:
    $$2(X/4)^{\frac{1}{3}} \cdot 2(X/27)^{\frac{1}{2}} \cdot \frac{1}{q^2} \prod_{p \nmid q} \left(1 - \frac{1}{p^{10}} \right) \sum_{\substack{a, b \in \mathbb{Z}/q}}  a_n(E_{a, b})  $$
    $$+ O \left( \frac{1}{X^{\frac{3}{4}}} \right) \cdot 2(X/4)^{\frac{1}{3}} \cdot 2(X/27)^{\frac{1}{2}} \cdot \frac{1}{q^2} \sum_{\substack{a, b \in \mathbb{Z}/q}}  a_n(E_{a, b}).$$
    The first of these terms is precisely the main term we seek to prove, while the second term is at most $O_\epsilon\left( \frac{1}{X^{\frac{3}{4}}} q^{- \frac{1}{2} + \epsilon}\prod_{p : p^2 \mid n} p^{\frac{v_p(n)}{2} + \epsilon} \right)$ 
    by Lemma \ref{Improved Trace Bound}, in the case where $y = z = 0.$ This term can later be subsumed into the main error term.

    Meanwhile, the same proof (replacing powers of $10$ by powers of $6$) shows that
    $$2(X/27)^{\frac{1}{2}} \cdot \frac{1}{q^2} \sum_{y, z \in \mathbb{Z}/q} \left( \sum_{a, b \in \mathbb{Z}/q}  a_n(E_{a, b}) \cdot e\left( \frac{a y + b z}{q}  \right)\right)\left( \sum_{\substack{m \leq X^{\frac{1}{12}}\\ q \mid m^4 y\\ q \mid m^6 z}} \frac{|\mu(m)|}{m^{6}} \right)$$
    $$= \left( 2(X/27)^{\frac{1}{2}} \cdot \frac{1}{q^2} \prod_{p \nmid q} \left(1 + \frac{1}{p^{6}} \right) \sum_{\substack{a, b \in \mathbb{Z}/q}}  a_n(E_{a, b}) \right) + O \left( \frac{1}{X^{\frac{5}{12}}} \right) \cdot 2(X/27)^{\frac{1}{2}} \cdot q^{- \frac{1}{2} + \epsilon}\prod_{p : p^2 \mid n} p^{\frac{v_p(n)}{2} + \epsilon}.$$
    We are obliged to take the absolute value of the $\mu(m)$ term since $O\left( \frac{2(X/27)^{\frac{1}{2}}}{m^6} \right)$ is an overestimate of the actual error. Replacing $\mu(m)$ with $|\mu(m)|$ has the effect of replacing the product $\prod_{p \nmid m'} \left(1 - \frac{1}{p^{6}} \right)$ with $\prod_{p \nmid m'} \left(1 + \frac{1}{p^{6}} \right)$, as written above. However, we note that
    \begin{align*}
        \prod_{p \nmid q} \left(1 + \frac{1}{p^{6}} \right) &\leq \prod_{p \nmid q} \left(1 + \frac{1}{p^{6}} + \frac{1}{p^{12}} + \cdots\right)^{-1} \leq \prod_{p} \left(1 + \frac{1}{p^{6}} + \frac{1}{p^{12}} + \cdots\right)^{-1} = \zeta(6) = O(1),
    \end{align*}
    so Lemma \ref{Improved Trace Bound} bounds the first error term by $O_\epsilon\left( X^{\frac{1}{2}} q^{- \frac{1}{2} + \epsilon}\prod_{p : p^2 \mid n} p^{\frac{v_p(n)}{2} + \epsilon} \right).$
    The second error term is bounded by $O \left( X^\frac{1}{6} q^{-\frac{1}{2} + \epsilon}\prod_{p : p^2 \mid n} p^{\frac{v_p(n)}{2} + \epsilon} \right),$ which is negligible. Analogously, we can bound the contribution from $O\left( \frac{2(X/4)^{\frac{1}{3}}}{m^4} \right)$ as
    \begin{align*}
        &\quad 2(X/4)^{\frac{1}{3}} \cdot \frac{1}{q^2} \sum_{y, z \in \mathbb{Z}/q} \left( \sum_{a, b \in \mathbb{Z}/q}  a_n(E_{a, b}) \cdot e\left( \frac{a y + b z}{q}  \right)\right) \sum_{\substack{m \leq X^{\frac{1}{12}}\\ q \mid (m^4 y, m^6 z)}} \frac{|\mu(m)|}{m^{4}}\\
        &\leq \left( 2(X/4)^{\frac{1}{3}} \cdot \frac{1}{q^2} \prod_{p \nmid q} \left(1 + \frac{1}{p^{4}} \right) \sum_{\substack{a, b \in \mathbb{Z}/q}}  a_n(E_{a, b}) \right) + O \left( \frac{1}{X^{\frac{1}{4}}} \right) \cdot O_\epsilon\left( X^{\frac{1}{3}} q^{- \frac{1}{2} + \epsilon}\prod_{p : p^2 \mid n} p^{\frac{v_p(n)}{2} + \epsilon} \right)\\
        &\leq O_\epsilon\left( X^{\frac{1}{3}} q^{- \frac{1}{2} + \epsilon}\prod_{p : p^2 \mid n} p^{\frac{v_p(n)}{2} + \epsilon} \right) + O_\epsilon\left( X^{\frac{1}{12}} q^{- \frac{1}{2} + \epsilon}\prod_{p : p^2 \mid n} p^{\frac{v_p(n)}{2} + \epsilon} \right),
    \end{align*}
    both of which are negligible.
\end{proof}

The next three lemmata bound the error terms of the sum in Proposition \ref{sums of coefficients}.

\begin{lemma}\label{Half error term 4}
Let $q = q(n).$ Then, as $X$ and $n$ grow,
\begin{align*}
    &\quad \frac{1}{q^2} \sum_{y, z \in \mathbb{Z}/q} \left( \sum_{a, b \in \mathbb{Z}/q}  a_n(E_{a, b}) \cdot e\left( \frac{a y + b z}{q}  \right)\right) \sum_{\substack{m \leq X^{\frac{1}{12}}\\ q \mid m^4 y\\ q \nmid m^6 z}} \mu(m) \sum_{\substack{A, B \in \mathbb{Z} \\ \max(4|A|^3, 27|B|^2) \leq X \\ m^4 \mid A, m^6 \mid B}}  e \left( \frac{-A y - Bz}{q}  \right)\\
    &\leq O_\epsilon\left( X^{\frac{1}{3}} q^{\frac{1}{2} + \epsilon}\prod_{p : p^2 \mid n} p^{\frac{v_p(n)}{2} + \epsilon}\right)
\end{align*}
\end{lemma}

\begin{proof}
    We first apply the triangle inequality and Lemma \ref{Improved Trace Bound} to obtain
    \begin{align*}
        \frac{1}{q^2} \sum_{y, z \in \mathbb{Z}/q} &\left| \sum_{a, b \in \mathbb{Z}/q}  a_n(E_{a, b}) \cdot e\left( \frac{a y + b z}{q}  \right)\right| \left( \sum_{\substack{m \leq X^{\frac{1}{12}}\\ q \mid m^4 y\\ q \nmid m^6 z}} |\mu(m)| \cdot \left| \sum_{\substack{A, B \in \mathbb{Z} \\ \max(4|A|^3, 27|B|^2) \leq X \\ m^4 \mid A, m^6 \mid B}} \! \! \! \! \! \! \! \! \! \! \! \! \! \! \! e \left( \frac{-A y - Bz}{q}  \right) \right| \right)\\
        &\ll_\epsilon q^{\frac{3}{2} + \epsilon}\prod_{p : p^2 \mid n} p^{\frac{v_p(n)}{2} + \epsilon} \cdot \frac{1}{q^2} \sum_{y, z \in \mathbb{Z}/q} \sum_{\substack{m \leq X^{\frac{1}{12}}\\ q \mid m^4 y\\ q \nmid m^6 z\\ \text{$m$ square-free}}} \left| \sum_{\substack{A, B \in \mathbb{Z} \\ \max(4|A|^3, 27|B|^2) \leq X \\ m^4 \mid A, m^6 \mid B}}  e \left( \frac{-A y - Bz}{q}  \right) \right|.
    \end{align*}
    Focusing on the second summand, we get
    \begin{align*}
        \left| \sum_{\substack{A, B \in \mathbb{Z} \\ \max(4|A|^3, 27|B|^2) \leq X \\ m^4 \mid A, m^6 \mid B}}  \! \! \! \! \! \! \! \! \! \! \! \! \! \!e \left( \frac{-A y - Bz}{q}  \right) \right| &= \left| \sum_{\substack{A \in \mathbb{Z} \\ 4|A|^3 \leq X \\ m^4 \mid A}} \! \! \! \!  e \left( \frac{-A y}{q}  \right) \! \! \! \! \! \sum_{\substack{B \in \mathbb{Z} \\ 27|B|^2 \leq X \\ m^6 \mid B}}  \! \! \! \! \!e \left( \frac{- Bz}{q}  \right) \right| \leq \frac{2\left( \frac{2 (X/4)^{\frac{1}{3}}}{m^4}  + 2 \right)}{1 - e\left( \frac{-m^6 z}{q}  \right) },
    \end{align*}
    since the first sum is a summation over $1,$ and the second sum is a geometric series with first and last terms having absolute value $1$ and common ratio $e\left( \frac{-m^6 z}{q} \right) \neq 1.$ The Taylor series expansion $1 - e(X) = - i2\pi X + \frac{(2\pi X)^2}{2!} + i \frac{(2\pi X)^3}{3!} - \frac{(2\pi X)^4}{4!} - \cdots$ implies that $|1 - e(X)| \geq \frac{1}{K} |X|$ over the interval $\left[ -\frac{1}{2}, \frac{1}{2} \right]$ for some sufficiently large $K.$ This implies that over $\left[ -\frac{1}{2}, \frac{1}{2} \right],$ we have $\left| \frac{1}{1 - e(X)}  \right| \leq \left| \frac{K}{X}  \right|.$ This allows us to bound the entire expression of Lemma \ref{Half error term 4} by
    $$q^{\frac{3}{2} + \epsilon}\prod_{p : p^2 \mid n} p^{\frac{v_p(n)}{2} + \epsilon} \cdot \frac{1}{q^2} \sum_{y, z \in \mathbb{Z}/q} \sum_{\substack{m \leq X^{\frac{1}{12}}\\ q \mid m^4 y\\ q \nmid m^6 z\\ \text{$m$ square-free}}} \left( \frac{2 (X/4)^{\frac{1}{3}}}{m^4}  + 2 \right)\left| \frac{2K q}{\overline{-m^6 z} } \right|,$$
    where $\overline{-m^6 z}$ denotes the reduction of $-m^6 z$ $\pmod q$ in the interval $\left[ - \frac{q}{2}, \frac{q}{2} \right).$ Switching the order of summations and rewriting each $m$ as $m' k$ as in the proof of Lemma \ref{Main Term}, with $\gcd (m', q) = 1$ and $k \mid q,$ we obtain
    $$2K q^{\frac{1}{2} + \epsilon}\prod_{p : p^2 \mid n} p^{\frac{v_p(n)}{2} + \epsilon} \sum_{\substack{k \mid q\\ \text{$k$ square-free}}} \frac{1}{k^4} \sum_{\substack{m' \leq X^{\frac{1}{12}}/k \\ \text{$m'$ square-free}\\ \gcd (m', q) = 1}} \left( \frac{2 (X/4)^{\frac{1}{3}}}{(m')^4}  + 2 \right) \sum_{\substack{z \in \mathbb{Z}/q \\ q \nmid k^6 z}} \left| \frac{1}{\overline{-(m' k)^6 z}}  \right| \sum_{\substack{y \in \mathbb{Z}/q \\ q \mid k^4 y}} 1.$$
    Because $m'$ is a unit in $\mathbb{Z}/q,$ we have
    \begin{align*}
        \sum_{\substack{z \in \mathbb{Z}/q \\ q \nmid k^6 z}} \left| \frac{1}{\overline{-(m' k)^6 z}} \right| &= \gcd(k^6, q) \sum_{0 < |z'| < q/2(k^6, q)} \frac{1}{|z'| \gcd (k^6, q)} \\
        &\leq 2\log\left( \frac{q}{2 \gcd (k^6, q)}  \right) + 1 \leq 2\log\left( q\right) + 1.
    \end{align*}
    Since we can factor out $2\log(q) + 1$ from the sum and lump it into the $q^\epsilon$ factor, our main term is of the form
    \begin{align*}
        &2K q^{\frac{1}{2} + \epsilon}\prod_{p : p^2 \mid n} p^{\frac{v_p(n)}{2} + \epsilon} \sum_{\substack{k \mid q\\ \text{$k$ square-free}}} \frac{1}{k^4} \sum_{\substack{m' \leq X^{\frac{1}{12}}/k \\ \text{$m'$ square-free} \\ \gcd(m', q) = 1}} \left( \frac{2 (X/4)^{\frac{1}{3}}}{(m')^4}  + 2 \right) \sum_{\substack{y \in \mathbb{Z}/q \\ q \mid k^4 y}} 1 \\
        \leq &2K q^{\frac{1}{2} + \epsilon}\prod_{p : p^2 \mid n} p^{\frac{v_p(n)}{2} + \epsilon} \sum_{\substack{k \mid q\\ \text{$k$ square-free}}} \frac{1}{k^4} \left( \sum_{\substack{m' \geq 1 \\ \text{$m'$ square-free} \\ \gcd (m', q) = 1}} \frac{2 (X/4)^{\frac{1}{3}}}{(m')^4}  + \sum_{\substack{m' \leq X^{\frac{1}{12}}/k}} 2 \right) \sum_{\substack{y \in \mathbb{Z}/q \\ q \mid k^4 y}} 1\\
        \leq &2K q^{\frac{1}{2} + \epsilon}\prod_{p : p^2 \mid n} p^{\frac{v_p(n)}{2} + \epsilon} \sum_{\substack{k \mid q\\ \text{$k$ square-free}}} \frac{1}{k^4} \left( 2 (X/4)^{\frac{1}{3}} \prod_{p \nmid q} \left(1 + \frac{1}{p^4} \right)  + \frac{2 X^{\frac{1}{12}}}{k} \right) \sum_{\substack{y \in \mathbb{Z}/q \\ q \mid k^4 y}} 1.
    \end{align*}
    Recall from the proof of Lemma \ref{Main Term} that $q \mid k^4 y$ imposes the same condition on $y$ as $c_y(k) \mid y.$ Because of this, we can rewrite 
    $$\sum_{\substack{y \in \mathbb{Z}/q \\ q \mid k^4 y}} 1 = \left( \prod_{\substack{p \mid k \\ p \in \{2, 3\}}} p^4  \right) \prod_{\substack{p \mid k\\ p \notin \{2, 3\}}} p.$$
    This allows us to rewrite the entire sum as
    \begin{align*}
        &2K q^{\frac{1}{2} + \epsilon}\prod_{p : p^2 \mid n} p^{\frac{v_p(n)}{2} + \epsilon} \sum_{\substack{k \mid q\\ \text{$k$ square-free}}} \frac{1}{k^4} \left( 2 (X/4)^{\frac{1}{3}} \prod_{p \nmid q} \left(1 + \frac{1}{p^4} \right)  + \frac{2 X^{\frac{1}{12}}}{k} \right) \left( \prod_{\substack{p \mid k \\ p \in \{2, 3\}}} p^4  \right) \prod_{\substack{p \mid k\\ p \notin \{2, 3\}}} p\\
        = & 4K q^{\frac{1}{2} + \epsilon}\prod_{p : p^2 \mid n} p^{\frac{v_p(n)}{2} + \epsilon} (X/4)^{\frac{1}{3}} \prod_{\substack{p \mid q\\ p \in \{2, 3\}}} (1 + p)\prod_{\substack{p \mid q\\ p \notin \{2, 3\}}} \left( 1 + \frac{1}{p^3}  \right) \prod_{p \nmid q} \left(1 + \frac{1}{p^4} \right)\\
        + & 4K q^{\frac{1}{2} + \epsilon}\prod_{p : p^2 \mid n} p^{\frac{v_p(n)}{2} + \epsilon} X^{\frac{1}{12}} \prod_{\substack{p \mid q\\ p \in \{2, 3\}}} (1 + p^0)\prod_{\substack{p \mid q\\ p \notin \{2, 3\}}} \left( 1 + \frac{1}{p^4}  \right).
    \end{align*}
    All the products over $p \mid q$ and $p \nmid q$ converge to constants. Therefore, the total sum is of the form $O_\epsilon\left(X^{\frac{1}{3}} q^{\frac{1}{2} + \epsilon}\prod_{p : p^2 \mid n} p^{\frac{v_p(n)}{2} + \epsilon} \right) +  O_\epsilon \left( X^{\frac{1}{12}} q^{\frac{1}{2} + \epsilon}\prod_{p : p^2 \mid n} p^{\frac{v_p(n)}{2} + \epsilon} \right).$ The second term is negligible.
\end{proof}

\begin{lemma}\label{Half error term 6}
    Let $q = q(n).$ Then, as $X$ and $n$ grow,
    \begin{align*}
        & \quad \frac{1}{q^2} \sum_{y, z \in \mathbb{Z}/q} \left( \sum_{a, b \in \mathbb{Z}/q}  a_n(E_{a, b}) \cdot e\left( \frac{a y + b z}{q}  \right)\right) \! \!\sum_{\substack{m \leq X^{\frac{1}{12}}\\ q \nmid m^4 y\\ q \mid m^6 z}} \mu(m) \! \! \! \! \! \! \! \! \! \! \!\sum_{\substack{A, B \in \mathbb{Z} \\ \max(4|A|^3, 27|B|^2) \leq X \\ m^4 \mid A, m^6 \mid B}} \! \! \! \! \! \! \! \! \! \! \! \! \! e \left( \frac{-A y - Bz}{q}  \right) \\
        &\leq O_\epsilon\left(X^{\frac{1}{2}} q^{\frac{1}{2} + \epsilon}\prod_{p : p^2 \mid n} p^{\frac{v_p(n)}{2} + \epsilon} \right) .
    \end{align*}
\end{lemma}

\begin{proof}
    The method is the same as for Lemma \ref{Half error term 4}. Again, we start by applying the triangle inequality and Lemma \ref{Improved Trace Bound} to get that the sum is
    \begin{align*}
        &\ll_\epsilon q^{\frac{3}{2} + \epsilon}\prod_{p : p^2 \mid n} p^{\frac{v_p(n)}{2} + \epsilon} \cdot \frac{1}{q^2} \sum_{y, z \in \mathbb{Z}/q} \sum_{\substack{m \leq X^{\frac{1}{12}}\\ q \nmid m^4 y\\ q \mid m^6 z\\ \text{$m$ square-free}}} \left| \sum_{\substack{A, B \in \mathbb{Z} \\ \max(4|A|^3, 27|B|^2) \leq X \\ m^4 \mid A, m^6 \mid B}}  e \left( \frac{-A y - Bz}{q}  \right) \right| \\
        &\leq q^{\frac{3}{2} + \epsilon}\prod_{p : p^2 \mid n} p^{\frac{v_p(n)}{2} + \epsilon} \cdot \frac{1}{q^2} \sum_{y, z \in \mathbb{Z}/q} \sum_{\substack{m \leq X^{\frac{1}{12}}\\ q \nmid m^4 y\\ q \mid m^6 z\\ \text{$m$ square-free}}} \left( \frac{2 (X/27)^{\frac{1}{2}}}{m^6}  + 2 \right) \cdot \frac{2}{1 - e\left( \frac{-m^4 y}{q}  \right) } \\
        &\leq q^{\frac{3}{2} + \epsilon}\prod_{p : p^2 \mid n} p^{\frac{v_p(n)}{2} + \epsilon} \cdot \frac{1}{q^2} \sum_{y, z \in \mathbb{Z}/q} \sum_{\substack{m \leq X^{\frac{1}{12}}\\ q \nmid m^4 y\\ q \mid m^6 z\\ \text{$m$ square-free}}} \left( \frac{2 (X/27)^{\frac{1}{2}}}{m^6}  + 2 \right) \cdot \left| \frac{2K q}{\overline{-m^4 y} } \right|,\\
        & \qquad \text{with $K$ as defined before, and $\overline{-m^4 z}$ being the reduction to the interval $\left[ -\frac{q}{2}, \frac{q}{2} \right)$}\\
        &= 2K q^{\frac{1}{2} + \epsilon}\prod_{p : p^2 \mid n} p^{\frac{v_p(n)}{2} + \epsilon} \sum_{\substack{k \mid q\\ \text{$k$ square-free}}} \frac{1}{k^6} \sum_{\substack{m' \leq X^{\frac{1}{12}}/k \\ \text{$m'$ square-free}\\ \gcd(m', q) = 1}} \left( \frac{2 (X/27)^{\frac{1}{2}}}{(m')^6}  + 2 \right) \sum_{\substack{z \in \mathbb{Z}/q \\ q \nmid k^4 y}} \left| \frac{1}{\overline{-(m' k)^4 y}}  \right| \sum_{\substack{z \in \mathbb{Z}/q \\ q \mid k^6 z}} 1\\
        & \qquad \text{after factoring $m$ as $m' \cdot k$ with $k$ dividing $q$ and $\gcd(m', k) = 1$}\\
        &\ll_\epsilon 2K q^{\frac{1}{2} + \epsilon}\prod_{p : p^2 \mid n} p^{\frac{v_p(n)}{2} + \epsilon} \sum_{\substack{k \mid q\\ \text{$k$ square-free}}} \frac{1}{k^6} \sum_{\substack{m' \leq X^{\frac{1}{12}}/k \\ \text{$m'$ square-free}\\ (m', q) = 1}} \left( \frac{2 (X/27)^{\frac{1}{2}}}{(m')^6}  + 2 \right) \sum_{\substack{z \in \mathbb{Z}/q \\ q \mid k^6 z}} 1\\
        & \qquad \text{after bounding and factoring out $\sum_{\substack{z \in \mathbb{Z}/q \\ q \nmid k^4 y}} \left| \frac{1}{\overline{-(m' y)^4 z}}  \right|  \leq 2 \log(q) + 1 \ll_\epsilon q^\epsilon$}\\
        &\leq 2K q^{\frac{1}{2} + \epsilon}\! \!  \prod_{p : p^2 \mid n} p^{\frac{v_p(n)}{2} + \epsilon} \! \! \! \! \! \sum_{\substack{k \mid q\\ \text{$k$ square-free}}} \! \! \! \! \! \frac{1}{k^6} \left( 2(X/27)^{\frac{1}{2}} \prod_{p \nmid q} \left( 1 + \frac{1}{p^6} \right) + \frac{2X^{\frac{1}{12}}}{k} \right) \cdot \left( \prod_{\substack{p \mid k \\ p \in \{2, 3\}}} p^6  \right)\prod_{\substack{p \mid k\\ p \notin \{2, 3\}}} p\\
        & \qquad \text{using the fact that the condition $q \mid k^6 z$ is equivalent to $c_z(k) \mid z$}\\
        &\ll_\epsilon 4K q^{\frac{1}{2} + \epsilon} \prod_{p : p^2 \mid n}  p^{\frac{v_p(n)}{2} + \epsilon} (X/27)^{\frac{1}{2}}  \prod_{\substack{p \mid q\\ p \in \{2, 3\}}} \left(1 + p^0 \right) \prod_{\substack{p \mid q\\ p \notin \{2, 3\}}} \left(1 + \frac{1}{p^5} \right) \!\prod_{p \nmid q} \left( 1 + \frac{1}{p^6} \right)\\
        &\ll_\epsilon X^{\frac{1}{2}} q^{\frac{1}{2} + \epsilon}\prod_{p : p^2 \mid n} p^{\frac{v_p(n)}{2} + \epsilon}
    \end{align*}
    We get the penultimate step by ignore the term $\frac{2X^{\frac{1}{12}}}{k},$ whose contribution is negligible, and factoring the remaining expression.
\end{proof}

\begin{lemma}\label{Log error term}
    Let $q = q(n).$ Then, as $X$ and $n$ grow,
    \begin{align*}
        & \quad \frac{1}{q^2} \sum_{y, z \in \mathbb{Z}/q} \left( \sum_{a, b \in \mathbb{Z}/q}  a_n(E_{a, b}) \cdot e\left( \frac{a y + b z}{q}  \right)\right) \! \!  \sum_{\substack{m \leq X^{\frac{1}{12}}\\ q \nmid m^4 y\\ q \nmid m^6 z}}  \! \!\mu(m) \! \! \! \! \! \! \! \! \! \! \! \! \sum_{\substack{A, B \in \mathbb{Z} \\ \max(4|A|^3, 27|B|^2) \leq X \\ m^4 \mid A, m^6 \mid B}}   \! \! \! \! \! \! \! \! \! \! \! \! \! \! \!e \left( \frac{-A y - Bz}{q}  \right)  \\
        & \leq O_\epsilon \left( X^{\frac{1}{12}} q^{\frac{3}{2} + \epsilon}\prod_{p : p^2 \mid n} p^{\frac{v_p(n)}{2} + \epsilon} \right).
    \end{align*}
\end{lemma}

\begin{proof}
    We can factor the final sum as 
    \begin{align*}
        \left| \sum_{\substack{A, B \in \mathbb{Z} \\ \max(4|A|^3, 27|B|^2) \leq X \\ m^4 \mid A, m^6 \mid B}}  e \left( \frac{-A y - Bz}{q}  \right) \right| &= \left|  \sum_{\substack{A \in \mathbb{Z} \\ 4|A|^3 \leq X \\ m^4 \mid A}} e\left( \frac{-Ay}{q}  \right) \right| \left| \sum_{\substack{B \in \mathbb{Z} \\ 27|B|^2 \leq X \\ m^6 \mid B}} e\left( \frac{-Bz}{q}  \right) \right|\\
        &\leq \left| \frac{2}{1 - e\left( \frac{-m^4 y}{q}  \right)}  \right| \left| \frac{2}{1 - e\left( \frac{-m^6 z}{q}  \right)}  \right| \\
        &\leq  \frac{2Kq}{ \left| \overline{-m^4 y}  \right|} \cdot  \frac{2Kq}{\left| \overline{-m^6 z}  \right| },
    \end{align*}
    where $K$ and $\overline{-m^4 y}$ are defined as in Lemmata \ref{Half error term 4} and \ref{Half error term 6}. Bounding the sum over $a, b \in \mathbb{Z}/q$ by Lemma \ref{Improved Trace Bound} and rearranging orders of summation, the full terms are bounded in absolute value by
    \begin{align*}
        &\ll_\epsilon  4K^2 \cdot q^{\frac{3}{2} + \epsilon}\prod_{p : p^2 \mid n} p^{\frac{v_p(n)}{2} + \epsilon} \sum_{\substack{m \leq X^{\frac{1}{12}}\\  \text{$m$ square-free}}} \sum_{\substack{y \in \mathbb{Z}/q \\ q \nmid m^4 y}}\frac{1}{ \left| \overline{-m^4 y}  \right|} \sum_{\substack{z \in \mathbb{Z}/q\\ q \nmid m^6 z}}  \frac{1}{\left| \overline{-m^6 z}  \right| }\\
        &\leq  4K^2 \cdot q^{\frac{3}{2} + \epsilon}\prod_{p : p^2 \mid n} p^{\frac{v_p(n)}{2} + \epsilon} \sum_{\substack{k \mid q\\\text{$k$ square-free}}} \sum_{\substack{m' \leq X^{\frac{1}{12}}/k\\  \text{$m'$ square-free}\\ \gcd(m', q) = 1}} \sum_{\substack{y \in \mathbb{Z}/q \\ q \nmid k^4 y}}\frac{1}{ \left| \overline{-(m')^4 k^4 y}  \right|} \sum_{\substack{z \in \mathbb{Z}/q\\ q \nmid k^6 z}}  \frac{1}{\left| \overline{-(m')^6 k^6 z}  \right| }\\
        & \leq 4K^2 \cdot q^{\frac{3}{2} + \epsilon}\prod_{p : p^2 \mid n} p^{\frac{v_p(n)}{2} + \epsilon} \sum_{\substack{k \mid q\\\text{$k$ square-free}}} \sum_{\substack{m' \leq X^{\frac{1}{12}}/k\\  \text{$m'$ square-free}\\ \gcd(m', q) = 1}} 1\\
        & \qquad \text{by bounding $\sum_{\substack{y \in \mathbb{Z}/q \\ q \nmid k^4 y}}\frac{1}{ \left| \overline{-(m')^4 k^4 y}  \right|} \sum_{\substack{z \in \mathbb{Z}/q\\ q \nmid k^6 z}}  \frac{1}{\left| \overline{-(m')^6 k^6 z}  \right| } \leq \left( 2\log(q) + 1 \right)^2 \ll_\epsilon q^\epsilon$} \\
        &\leq 4K^2 \cdot q^{\frac{3}{2} + \epsilon}\prod_{p : p^2 \mid n} p^{\frac{v_p(n)}{2} + \epsilon} X^{\frac{1}{12}} \sum_{\substack{k \mid q\\\text{$k$ square-free}}} \frac{1}{k}\\
        &\ll_\epsilon  X^{\frac{1}{12}} q^{\frac{3}{2} + \epsilon}\prod_{p : p^2 \mid n} p^{\frac{v_p(n)}{2} + \epsilon},
    \end{align*}
    by the chain of inequalities
    \begin{align*}
        \sum_{\substack{k \mid q\\\text{$k$ square-free}}} \frac{1}{k} &= \prod_{p \mid q} \left(1 + \frac{1}{p} \right) \leq \prod_{p \mid q} \left(1 - \frac{1}{p} \right)^{-1} = \frac{q}{\phi(q)} \leq O\left( \log \log q  \right),
    \end{align*}
    the last inequality being due to a classical result \cite[Theorem 328]{HW79}.
\end{proof}

We observe that, together, the preceding lemmata imply Proposition \ref{sums of coefficients}.

\begin{proof}[Proof of Proposition \ref{sums of coefficients}]
    Lemma \ref{Polya Vinogradov} reformulates  $\sum_{\substack{A, B \in \mathbb{Z} \\ \max(4|A|^3, 27|B|^2) \leq X \\ \text{$p^4 \nmid A$ or $p^6 \nmid b$ for all $p$}\\ 4A^3 + 27 B^2 \neq 0}} a_n(E_{A, B}) $ as
    $$\frac{1}{q^2} \sum_{y, z \in \mathbb{Z}/q} \left( \sum_{a, b \in \mathbb{Z}/q}  a_n(E_{a, b}) \cdot e\left( \frac{a y + b z}{q}  \right)\right)\left( \sum_{\substack{A, B \in \mathbb{Z} \\ \max(4|A|^3, 27|B|^2) \leq X \\ \text{$p^4 \nmid A$ or $p^6 \nmid b$ for all $p$}}}  e \left( \frac{-A y - Bz}{q}  \right) \right),$$
    plus an error term of $O_\epsilon \left( X^{\frac{1}{6}} n^{\frac{1}{2} + \epsilon} \right).$ Using inclusion-exclusion, we can rewrite
    $$\sum_{\substack{A, B \in \mathbb{Z} \\ \max(4|A|^3, 27|B|^2) \leq X \\ \text{$p^4 \nmid A$ or $p^6 \nmid b$ for all $p$}}}  e \left( \frac{-A y - Bz}{q}  \right) = \sum_{m \leq X^{\frac{1}{12}}} \mu(m) \sum_{\substack{A, B \in \mathbb{Z} \\ \max(4|A|^3, 27|B|^2) \leq X \\ m^4 \mid A, m^6 \mid B}}  e \left( \frac{-A y - Bz}{q}  \right).$$
    We split the main term into four sums over different cases of $m$: (1) $q \mid m^4 y$ and $q \mid m^6 z,$ (2) $q \mid m^4y$ but $q \nmid m^6 z$, (3) $q \mid m^6 z$ but $q \nmid m^4 y,$ or (4) $q \nmid m^4 y$ and $q \nmid m^6 z.$ Lemmata \ref{Main Term}, \ref{Half error term 4}, \ref{Half error term 6}, and \ref{Log error term} together give that this sum is equal to
    $$\frac{1}{q^2} \cdot 2(X/4)^{\frac{1}{3}} \cdot 2(X/27)^{\frac{1}{2}} \cdot \prod_{p \nmid q} \left( 1 - \frac{1}{p^{10}} \right) \cdot \sum_{a, b \in \mathbb{Z}/q} a_n(E_{a, b}) + O_\epsilon\left( X^{\frac{1}{2}} q^{- \frac{1}{2} + \epsilon}\prod_{p : p^2 \mid n} p^{\frac{v_p(n)}{2} + \epsilon} \right)$$
    $$+ O_\epsilon\left( X^{\frac{1}{3}} q^{\frac{1}{2} + \epsilon}\prod_{p : p^2 \mid n} p^{\frac{v_p(n)}{2} + \epsilon}\right) + O_\epsilon\left(X^{\frac{1}{2}} q^{\frac{1}{2} + \epsilon}\prod_{p : p^2 \mid n} p^{\frac{v_p(n)}{2} + \epsilon} \right) + O_\epsilon \left( X^{\frac{1}{12}} q^{\frac{3}{2} + \epsilon}\prod_{p : p^2 \mid n} p^{\frac{v_p(n)}{2} + \epsilon} \right).$$
    The error terms $O_\epsilon\left( X^{\frac{1}{2}} q^{- \frac{1}{2} + \epsilon}\prod_{p : p^2 \mid n} p^{\frac{v_p(n)}{2} + \epsilon} \right)$ and $ O_\epsilon\left( X^{\frac{1}{3}} q^{\frac{1}{2} + \epsilon}\prod_{p : p^2 \mid n} p^{\frac{v_p(n)}{2} + \epsilon}\right)$ are subsumed by the term $O_\epsilon\left(X^{\frac{1}{2}} q^{\frac{1}{2} + \epsilon}\prod_{p : p^2 \mid n} p^{\frac{v_p(n)}{2} + \epsilon} \right)$ and can be removed. 
\end{proof}

\subsection{The Proof of Theorem \ref{Main Result}}

We begin with the following lemma, which combined with Proposition \ref{sums of coefficients} will help us prove Lemma \ref{Average of coefficients}.

\begin{lemma}\label{Number of curve of bounded height}
    As $X$ grows,
    $$\sum_{\substack{A, B \in \mathbb{Z} \\ \max(4|A|^3, 27|B|^2) \leq X \\ \text{$p^4 \nmid A$ or $p^6 \nmid b$ for all $p$}\\ 4A^3 + 27 B^2 \neq 0}} 1 = 2(X/4)^{\frac{1}{3}} \cdot 2(X/27)^{\frac{1}{2}} \prod_{p} \left( 1 - \frac{1}{p^{10}} \right) +  O\left( X^{\frac{1}{2}}  \right).$$
\end{lemma}

\begin{proof}
    This is a simplified version of \cite[Lemma 13]{SS25}. We can drop the condition of $4A^3 + 27 B^2 \neq 0,$ which introduces a negligible error of $O \left(X^{\frac{1}{6}} \right).$ Inclusion-exclusion gives
    \begin{align*}
        \sum_{\substack{A, B \in \mathbb{Z} \\ \max(4|A|^3, 27|B|^2) \leq X \\ \text{$p^4 \nmid A$ or $p^6 \nmid b$ for all $p$}\\ 4A^3 + 27 B^2 \neq 0}} 1 &= \sum_{m < X^{\frac{1}{12}}} \mu(m) \sum_{\substack{A, B \in \mathbb{Z} \\ \max(4|A|^3, 27|B|^2) \leq X \\ m^4 \mid A\\ m^6 \mid B}} 1 \\
        &= \sum_{m < X^{\frac{1}{12}}} \mu(m) \left( \frac{2(X/4)^{\frac{1}{3}} \cdot 2(X/27)^{\frac{1}{2}}}{m^{10}}  + O\left( \frac{X^{\frac{1}{2}}}{m^6}  \right) + O\left( \frac{X^{\frac{1}{3}}}{m^4}  \right)\right) \\
        &= 2(X/4)^{\frac{1}{3}} \cdot 2(X/27)^{\frac{1}{2}} \sum_{m < X^{\frac{1}{12}}} \frac{\mu(m)}{m^{10}} + O\left( X^{\frac{1}{2}}  \right)\\
        &= 2(X/4)^{\frac{1}{3}} \cdot 2(X/27)^{\frac{1}{2}} \sum_{m} \frac{\mu(m)}{m^{10}} + O\left( \frac{1}{X^{\frac{3}{4}}}  \right) + O\left( X^{\frac{1}{2}}  \right)\\
        &= 2(X/4)^{\frac{1}{3}} \cdot 2(X/27)^{\frac{1}{2}} \prod_{p} \left( 1 - \frac{1}{p^{10}} \right) +  O\left( X^{\frac{1}{2}}  \right). \qedhere
    \end{align*}

\begin{proof}[Proof of Lemma \ref{Average of coefficients}]
    Given two functions $F$ and $G$ and with error terms $\epsilon_1$ and $\epsilon_2$ respectively, where $\epsilon_2 \ll G,$ we can use a geometric series to calculate the ratio
    $$\frac{F + \epsilon_1}{G + \epsilon_2} = \frac{F + \epsilon_1}{G} \cdot \frac{1}{1 + \frac{\epsilon_2}{G}} = \frac{F + \epsilon_1}{G} \left( 1 - \frac{\epsilon_2}{G} + \frac{\epsilon_2^2}{G^2} - \cdots \right) = \frac{F}{G} + O\left( \frac{\epsilon_1}{G} \right) + O\left( \frac{F \epsilon_2}{G^2} \right).$$
    We apply Proposition \ref{sums of coefficients} and Lemma \ref{Number of curve of bounded height} to set
    \begin{align*}
        F &= \frac{1}{q^2} \cdot 2(X/4)^{\frac{1}{3}} \cdot 2(X/27)^{\frac{1}{2}} \cdot \prod_{p \nmid q} \left( 1 - \frac{1}{p^{10}} \right) \cdot \sum_{a, b \in \mathbb{Z}/q} a_n(E_{a, b})\\
        \epsilon_1 &\leq O_\epsilon \left( X^{\frac{1}{6}} n^{\frac{1}{2} + \epsilon} \right) + O_\epsilon\left(X^{\frac{1}{2}} q^{\frac{1}{2} + \epsilon}\prod_{p : p^2 \mid n} p^{\frac{v_p(n)}{2} + \epsilon} \right) + O_\epsilon \left( X^{\frac{1}{12}} q^{\frac{3}{2} + \epsilon}\prod_{p : p^2 \mid n} p^{\frac{v_p(n)}{2} + \epsilon} \right)\\
        G &= 2(X/4)^{\frac{1}{3}} \cdot 2(X/27)^{\frac{1}{2}} \prod_{p} \left( 1 - \frac{1}{p^{10}} \right)\\
        \epsilon_2 &\leq O\left( X^{\frac{1}{2}}  \right).
    \end{align*}
    Therefore,
    \begin{align*}
        \mathbb{E}_{\{E: H(E) \leq X\}}[a_n(E)] &= \mathbb{E}_{\substack{A, B \in \mathbb{Z} \\ \max(4|A|^3, 27|B|^2) \leq X \\ \text{$p^4 \nmid A$ or $p^6 \nmid b$ for all $p$}\\ 4A^3 + 27 B^2 \neq 0}} [a_n(E)]\\
        &= \frac{1}{q^2} \prod_{p \mid q} \left( 1 - \frac{1}{p^{10}} \right)^{-1} \cdot \sum_{a, b } a_n(E_{a, b}) + O \left( X^{-\frac{2}{3}} n^{\frac{1}{2}} \right)\\
        &+ O_\epsilon\left(X^{-\frac{1}{3}} q^{\frac{1}{2} + \epsilon}\prod_{p : p^2 \mid n} p^{\frac{v_p(n)}{2} + \epsilon} \right) + O_\epsilon \left( X^{-\frac{3}{4}} q^{\frac{3}{2} + \epsilon}\prod_{p : p^2 \mid n} p^{\frac{v_p(n)}{2} + \epsilon} \right)\\
        &+ O \left( X^{-\frac{1}{3}} \cdot \frac{1}{q^2} \cdot \sum_{a, b \in \mathbb{Z}/q} a_n(E_{a, b}) \right).
    \end{align*}
    Lemma \ref{Improved Trace Bound} for $y = z = 0$ shows $ X^{-\frac{1}{3}} \cdot \frac{1}{q^2} \cdot \sum_{a, b \in \mathbb{Z}/q} a_n(E_{a, b}) \ll_\epsilon X^{-\frac{1}{3}} q^{-\frac{1}{2} + \epsilon}\prod_{p : p^2 \mid n} p^{\frac{v_p(n)}{2} + \epsilon},$ so the last error term can be ignored. Finally, using the Chinese remainder theorem and the multiplicativity of $a_n,$ we write the main term as
    \begin{align*}
        \frac{1}{q^2} \cdot \prod_{p \mid q} \left( \frac{1}{1 - p^{-10}} \right) \cdot \sum_{a, b \in \mathbb{Z}/q} a_n(E_{a, b}) &= \prod_{p \mid q} \frac{1}{p^{2v_p(q)}} \cdot \frac{1}{1 - p^{-10}} \sum_{a, b \in \mathbb{Z}/p^{v_p(q)}} a_{p^{v_p(n)}}(E_{a, b})\\
        &= \prod_{p \mid q} \frac{1}{1 - p^{-10}} \int_{\mathbb{Z}_p^2 \setminus(p^4 \mathbb{Z}_p \times p^6 \mathbb{Z}_p)} a_{p^{v_p(n)}}(E_{a, b}) = \prod_{p \mid q} \ell_{p, v_p(n)}.
    \end{align*}
    Here, we are able to rewrite the sum as an integral over the uniform measure on $\mathbb{Z}_p^2$ with total mass $1$ in accordance with Remark \ref{p adic coefficients}.
\end{proof}

Finally, we use the preceding results to prove Theorem \ref{Main Result}.

\begin{proof}[Proof of Theorem \ref{Main Result}]
    We begin by following the proof of \cite[Proposition 7]{SS25}. The Voronoi summation formula (Lemma \ref{Voronoi Summation}) gives
    $$\frac{\epsilon(E)}{N(E)} \sum_{n = 1}^\infty a_n(E) W\left( L\left( \frac{n}{N(E)} - x\right)\right) = \sum_{n = 1}^{\infty} \frac{a_n(E)}{\sqrt{n}}  \int_0^\infty 2\pi\sqrt{u} W(L(u - x)) J_1\left( 4\pi \sqrt{un} \right) du.$$
    It then suffices to find the expectation of the right-hand side. Since the sum over $n$ converges absolutely as a consequence of \cite[Lemma 12]{SS25}, we can rearrange terms and get
    \begin{align}
        & \quad \sum_{n = 1}^{\infty} \frac{1}{\sqrt{n}}  \int_0^\infty 2\pi\sqrt{u} W(L(u - x)) J_1\left( 4\pi \sqrt{un} \right) du \cdot \mathbb{E}_{E : H(E) \leq X} \left[  a_n(E) \right] \nonumber\\
        &= \sum_{n = 1}^{\infty} \frac{1}{\sqrt{n}}  \int_0^\infty 2\pi\sqrt{u} W(L(u - x)) J_1\left( 4\pi \sqrt{un} \right) du \cdot \prod_{p \mid n} \ell_{p, v_p(n)} \nonumber\\
        &+ \sum_{n = 1}^{\infty} \frac{1}{\sqrt{n}}  \int_0^\infty 2\pi\sqrt{u} W(L(u - x)) J_1\left( 4\pi \sqrt{un} \right) du   \cdot  O_\epsilon\left(X^{-\frac{1}{3}} q^{\frac{1}{2} + \epsilon}\prod_{p : p^2 \mid n} p^{\frac{v_p(n)}{2} + \epsilon} \right) \label{eqn:error-term-1}\\ 
        &+ \sum_{n = 1}^{\infty} \frac{1}{\sqrt{n}}  \int_0^\infty 2\pi\sqrt{u} W(L(u - x)) J_1\left( 4\pi \sqrt{un} \right) du \cdot O_\epsilon \left( X^{-\frac{3}{4}} q^{\frac{3}{2} + \epsilon}\prod_{p : p^2 \mid n} p^{\frac{v_p(n)}{2} + \epsilon} \right) \label{eqn:error-term-2}\\
        &+ \sum_{n = 1}^{\infty} \frac{1}{\sqrt{n}}  \int_0^\infty 2\pi\sqrt{u} W(L(u - x)) J_1\left( 4\pi \sqrt{un} \right) du \cdot O \left( X^{-\frac{2}{3}} n^{\frac{1}{2}} \right). \label{eqn:error-term-3}
    \end{align}
    by Lemma \ref{Average of coefficients}. The end of the proof of \cite[Proposition 7]{SS25} establishes that the main term can be arranged to place the sum and the $\prod_{p \mid n} \ell_{p, v_p(n)}$ factor inside the integral. Meanwhile,
    \begin{align*}
        \eqref{eqn:error-term-1} &\leq \sum_{n = 1}^{\infty} O_{m, W}\left( L^{m - 1} n^{-\frac{m}{2} - \frac{1}{4}}  \right) \cdot \frac{1}{\sqrt{n}}  \cdot  O_\epsilon\left(X^{-\frac{1}{3}} q^{\frac{1}{2} + \epsilon}\prod_{p : p^2 \mid n} p^{\frac{v_p(n)}{2} + \epsilon} \right) \qquad \text{by Lemma \ref{Test function integral bound}}\\
        &\ll_{W, m, \epsilon} L^{m - 1} X^{-\frac{1}{3}} \sum_{n = 1}^\infty n^{-\frac{m}{2} - \frac{3}{4}} \cdot q^{\frac{1}{2} + \epsilon}\prod_{p : p^2 \mid n} p^{\frac{v_p(n)}{2} + \epsilon}\\
        &\ll L^{m - 1} X^{-\frac{1}{3}} \prod_{\text{$p$ prime}} \left( 1 + p^{- \frac{m}{2} - \frac{1}{4} + \epsilon} + \sum_{v \geq 2} p^{-\frac{v(2m + 1)}{4} + \frac{1}{2} + 2\epsilon} \right)\\
        & \qquad \text{(we have $\ll$ instead of $=$ because $v_2(q)$ and $v_3(q)$ can be as big as $7$)}\\
        &= L^{m - 1} X^{-\frac{1}{3}} \prod_{p} \left( 1 + p^{- \frac{m}{2} - \frac{1}{4} + \epsilon} + \frac{p^{-m + \frac{1}{4} + 2\epsilon}}{1 - p^{- \frac{2m + 1}{4}}} \right) \\
        & \leq L^{m - 1} X^{-\frac{1}{3}} \prod_{p} \left( 1 + p^{- \frac{m}{2} - \frac{1}{4} + \epsilon} + p^{-m + \frac{1}{4} + 2\epsilon} \right),
    \end{align*}
    which is $O_\epsilon \left(  L^{\frac{1}{2} + \epsilon} X^{-\frac{1}{3}} \right)$ since the product converges absolutely for $m > \frac{3}{2} + \epsilon$. Next,
    \begin{align*}
        \eqref{eqn:error-term-2} &\leq \sum_{n = 1}^{\infty} O_{m, W}\left( L^{m - 1} n^{- \frac{m}{2} - \frac{1}{4}}  \right) \cdot \frac{1}{\sqrt{n}}  \cdot  O_\epsilon \left( X^{-\frac{3}{4}} q^{\frac{3}{2} + \epsilon}\prod_{p : p^2 \mid n} p^{\frac{v_p(n)}{2} + \epsilon} \right)\\
        &\ll_{W, m, \epsilon} L^{m - 1} X^{-\frac{3}{4}} \sum_{n = 1}^\infty n^{- \frac{m}{2} - \frac{3}{4}} \cdot q^{\frac{3}{2} + \epsilon}\prod_{p : p^2 \mid n} p^{\frac{v_p(n)}{2} + \epsilon} \qquad \text{by Lemma \ref{Test function integral bound}}\\
        & \ll L^{m - 1}  X^{-\frac{3}{4}}\prod_{\text{$p$ prime}}  \left( 1 + p^{-\frac{m}{2} + \frac{3}{4} + \epsilon} + \sum_{v \geq 2} p^{-\frac{v(2m + 1)}{4} + \frac{3}{2} + 2\epsilon} \right)\\
        &= L^{m - 1}X^{-\frac{3}{4}} \prod_{p} \left( 1 + p^{-\frac{m}{2} + \frac{3}{4} + \epsilon} + \frac{p^{- m + \frac{5}{4} + 2\epsilon}}{1 - p^{- \frac{2m + 1}{4}}} \right) \\
        &\leq L^{m - 1}X^{-\frac{3}{4}}  \prod_{p} \left( 1 + p^{-\frac{m}{2} + \frac{3}{4} + \epsilon} + p^{- m + \frac{5}{4} + 2\epsilon} \right),
    \end{align*}
    which is $O_\epsilon \left( L^{- \frac{5}{2} + \epsilon} X^{-\frac{3}{4}} \right)$ since the product converges absolutely for $m > \frac{7}{2} + \epsilon.$ Finally,
    \begin{align*}
        \eqref{eqn:error-term-3} &\leq \sum_{n = 1}^{\infty} O_{m, W}\left( L^{m - 1} n^{- \frac{m}{2} - \frac{1}{4}}  \right) \cdot \frac{1}{\sqrt{n}}  \cdot  O \left( X^{-\frac{2}{3}} n^{\frac{1}{2}} \right) \qquad \text{by Lemma \ref{Test function integral bound}} \\
        &\ll L^{m - 1} X^{- \frac{2}{3}} \sum_{n = 1}^{\infty} O_{m, W}     \left( n^{- \frac{m - 1}{2} - \frac{3}{4}}  \right),
    \end{align*}
    which is $O_\epsilon \left( L^{\frac{1}{2} + \epsilon} X^{-\frac{2}{3}} \right)$ since the product converges absolutely for $m > \frac{3}{2}$. All three error terms thus decay when $L \ll X^{\frac{3}{10} - \epsilon}$ for some choice of $\epsilon.$
\end{proof}

\end{proof}

\bibliographystyle{amsplain}
\bibliography{refs}
    
\end{document}